%% file: dhn-adaptive.tex
\documentclass[a4paper,american,reqno]{amsart}

\newif\ifPreprint \Preprinttrue
\newif\ifSubmission \Submissionfalse

\input{setup-preprint}
\bibliography{dhn-adaptive}

\begin{document}

\title[An Adaptive Optimization Method for District Heating
Networks]{Adaptive Nonlinear Optimization of District Heating Networks
  Based on Model and Discretization Catalogs}

\author[H. Dänschel]{Hannes Dänschel}
\author[V. Mehrmann]{Volker Mehrmann}
\author[M. Roland]{Marius Roland}
\author[M. Schmidt]{Martin Schmidt}

\address[H. Dänschel, V. Mehrmann]
{TU Berlin, Inst.\ f.\ Mathematik, MA 4-5,
  Str.\ des 17.\ Juni 135, 10623 Berlin, Germany}

\email{daenschel@math.tu-berlin.de}
\email{mehrmann@math.tu-berlin.de}

\address[M. Roland, M. Schmidt]{Trier University,
  Department of Mathematics,
  Universitätsring 15,
  54296 Trier,
  Germany}

\email{roland@uni-trier.de}
\email{martin.schmidt@uni-trier.de}

\date{\today}

\begin{abstract}
  \input{abstract}
\end{abstract}

\keywords{\input{keywords}}
\subjclass[2010]{\input{msc2010}}

\maketitle

\input{introduction}
\input{modeling}
\input{error-measures}
\input{adaptive-algorithm}
\input{numerical-results}
\input{conclusion}
\input{acknowledgements}

\printbibliography

\end{document}

%%% Local Variables:
%%% mode: latex
%%% TeX-master: t
%%% End:

%% file: setup-preprint.tex
\usepackage{babel}
\usepackage[utf8]{inputenc}
\usepackage[T1]{fontenc}
\usepackage{lmodern}
\usepackage[binary-units=true]{siunitx}
\usepackage{booktabs}
\usepackage{xspace}
\usepackage{accents}
\usepackage[ruled,vlined,linesnumbered]{algorithm2e}
\usepackage{mathtools}
\usepackage{todonotes}
\usepackage{pgfplots}
\usepackage{tikz}
\usepackage{csquotes}
\usepackage{microtype}
\usepackage{eurosym}
\usepackage[style = numeric-comp,
            giveninits = true,
            maxbibnames = 100,
            maxcitenames = 2,
            isbn = false,
            backend = biber]{biblatex}
\usepackage[colorlinks,
            citecolor=blue,
            urlcolor=blue,
            linkcolor=blue]{hyperref} % should always be the last package

\makeatletter
\patchcmd{\@settitle}{\uppercasenonmath\@title}{\scshape\large}{}{}
\patchcmd{\@setauthors}{\MakeUppercase}{\scshape\normalsize}{}{}
\makeatother

\vfuzz \hfuzz
\usetikzlibrary{decorations.pathreplacing}
\tikzset{every picture/.style={remember picture}}
\pgfplotsset{compat=1.15}
\DeclareUnicodeCharacter{2212}{-}

\input{macros}

%% file: macros.tex
\newcommand{\JMR}[1]{\todo[author=MR,color=green!50,size=\small]{#1}}

\newcommand{\abbr}[1][abbrev]{#1.\xspace}

\newcommand{\cf}{\abbr[cf]}
\newcommand{\eg}{\abbr[e.g]}

\newcommand{\ie}{\abbr[i.e]}

\newcommand{\wrt}{\abbr[w.r.t]}

\newcommand{\field}{\mathbb}

\newcommand{\complex}{\field{C}}

\makeatletter
\@ifundefined{C}{\newcommand{\C}{\complex}}
{\PackageWarning{math-fields}{Use \string\field{C} for complex numbers;
 \MessageBreak command \string\C\space already defined}{}}
\makeatother

\newtheorem{assumption}{Assumption}
\newtheorem{lemma}{Lemma}
\newtheorem{remark}{Remark}
\newtheorem{theorem}{Theorem}
\newtheorem{definition}{Definition}

\newtheorem{corollary}{Corollary}

\providecommand{\sfnamesize}{\relsize{0}}
\providecommand{\sfname}[1]{\textsf{\sfnamesize#1}\xspace}

\newcommand{\CONOPTfour}{\sfname{CONOPT4}}

\newcommand{\GAMS}{\sfname{GAMS}}

\newcommand{\Python}{\sfname{Python}}
\newcommand{\Pyomo}{\sfname{Pyomo}}

\newcommand{\AROMA}{\sfname{AROMA}}
\newcommand{\STREET}{\sfname{STREET}}

\newcommand{\graphset}{}
\newcommand{\Arcs}{\graphset{A}}

\newcommand{\Nodes}{\graphset{V}}
\newcommand{\Vertices}{\Nodes}

\newcommand{\node}{u}
\newcommand{\otherNode}{v}

\newcommand{\arc}{a}

\newcommand{\Outedges}{\delta^{\text{out}}}

\newcommand{\Inedges}{\delta^{\text{in}}}

\newcommand{\solution}{y}

\newcommand{\Activeset}[1]{\mathcal{A}\ifx#1\empty\left(#1\right)\fi}

\newcommand{\activeset}[1]{\mathcal{A}\ifx#1\empty\else(#1)\fi}

\newcommand{\primalinfeasibility}[1][]{\theta_{\mathrm{pri}\ifx#1\empty\else{,{#1}}\fi}}

\newcommand{\red}{\reduced}
\newcommand{\kktrighthandside}[1][]{\omega\ifx#1\empty\else_{#1}\fi}
\newcommand{\kktrhs}{\kktrighthandside}
\newcommand{\kktreducedrighthandside}[1][\empty]{\kktrhs\ifx#1\empty_{\text{\red}}\else_{\text{\red,#1}}\fi}

\newcommand{\reduced}{\text{r}}

\newcommand{\steplength}[2][]{\alpha\ifx#1\empty_{#2}\else_{#2,#1}\fi}

\newcommand{\st}{\text{s.t.}}

\makeatletter
\newcommand{\fcdot}{\,\cdot\,}
\newcommand{\fcarg}[1]{\def\fc@rg{#1}\ifx\fc@rg\empty\fcdot\else\fc@rg\fi}
\makeatother
\newcommand{\abs}[1]{\lvert\fcarg{#1}\rvert}
\newcommand{\Abs}[1]{\left\lvert#1\right\rvert}

\newcommand{\Defset}[3][\defsep]{\Set{#2#1#3}}
\newcommand{\norm}[2][]{\lVert\fcarg{#2}\rVert\ifx#1\empty\else_{#1}\fi}
\newcommand{\Norm}[2][]{\left\lVert#2\right\rVert\ifx#1\empty\else_{#1}\fi}

\newcommand{\Set}[1]{\left\{#1\right\}}
\newcommand{\snorm}[2][]{\lvert\!\lvert\!\lvert
  \fcarg{#1}\rvert\!\rvert\!\rvert\ifx#2\empty\else_{#1}\fi}
\newcommand{\Snorm}[2][]{\left\lvert\!\left\lvert\!\left\lvert
  #2\right\rvert\!\right\rvert\!\right\rvert\ifx#1\empty\else_{#1}\fi}
\newcommand{\sprod}[3][]{%
  \langle\fcarg{#2},\fcarg{#3}\rangle\ifx#1\empty\else_{#1}\fi}
\newcommand{\Sprod}[3][]{%
  \left\langle\fcarg{#2},\fcarg{#3}\right\rangle\ifx#1\empty\else_{#1}\fi}

\newcommand{\optmathindex}[1]{\ifx#1\empty\else,#1\fi}
\newcommand{\opttextindex}[1]{\ifx#1\empty\else,\text{#1}\fi}
\newcommand{\optmathsb}[1]{\ifx#1\empty\else_{#1}\fi}
\newcommand{\opttextsb}[1]{\ifx#1\empty\else_{\text{#1}}\fi}
\newcommand{\optmathsp}[1]{\ifx#1\empty\else^{#1}\fi}
\newcommand{\opttextsp}[1]{\ifx#1\empty\else^{\text{#1}}\fi}

\newcommand{\continuousFunctions}[1]{\mathcal{C}\ifx#1\empty\else^{#1}\fi}

\newcommand{\piecewiseContinuousFunctions}[1]{\mathcal{C}_p\ifx#1\empty\else^{#1}\fi}

\newcommand{\define}{\mathrel{{\mathop:}{=}}}
\newcommand{\enifed}{\mathrel{{=}{\mathop:}}}

\newcommand{\diff}{\xspace\,\mathrm{d}}

\makeatletter
\newcommand{\objref}[4]{\def\obj@rg{#4}%
  #1\ifx\obj@rg\empty#2\else#3\xspace\ref{#4}--\fi\ref}
\makeatother

\newcommand{\Vdh}{\Nodes}

\newcommand{\Vff}{\Vdh_{\text{ff}}}

\newcommand{\Vbf}{\Vdh_{\text{bf}}}

\newcommand{\Aff}{\Arcs_{\text{ff}}}

\newcommand{\Abf}{\Arcs_{\text{bf}}}

\newcommand{\Apipe}{\Arcs_{\text{p}}}

\newcommand{\Acons}{\Arcs_{\text{c}}}

\newcommand{\nodeArcs}{\delta}

\newcommand{\cost}{C}
\newcommand{\wastecost}{\cost_{\text{w}}}
\newcommand{\gascost}{\cost_{\text{g}}}
\newcommand{\pressurecost}{\cost_{\text{p}}}

\newcommand{\densityinternalenergy}{e}
\newcommand{\densintene}{\densityinternalenergy}
\newcommand{\densityinternalentropy}{s}
\newcommand{\densintent}{\densityinternalentropy}

\newcommand{\heattransfercoefficient}{k_W}
\newcommand{\heattrans}{\heattransfercoefficient}

\newcommand{\negativemassflow}{\gamma}
\newcommand{\negmflow}{\negativemassflow}
\newcommand{\positivemassflow}{\beta}
\newcommand{\posmflow}{\positivemassflow}
\newcommand{\Powerwaste}{\Power_{\mathrm{w}}}
\newcommand{\Powergas}{\Power_{\mathrm{g}}}
\newcommand{\Powerpress}{\Power_{\mathrm{p}}}

\newcommand{\soiltemp}{\temp_\text{W}}

\newcommand{\stagnation}{\mathrm{s}}

\newcommand{\spacediff}{\Delta x}
\newcommand{\spacepointsset}{\varGamma}

\newcommand{\densintenebf}{\densintene^{\text{bf}}}

\newcommand{\densinteneforeflowarc}{\densintene^{\text{ff}}_\arc}

\newcommand{\densinteneffarc}{\densinteneforeflowarc}

\newcommand{\level}{\ell}
\newcommand{\new}{\text{new}}

\newcommand{\area}{A}
\newcommand{\bore}{D}% inner pipe diameter

\newcommand{\density}{\rho}
\newcommand{\dens}{\density}

\newcommand{\diameter}{\bore}
\newcommand{\diam}{\diameter}

\newcommand{\gravity}{g}
\newcommand{\grav}{\gravity}

\newcommand{\length}{L}
\newcommand{\massflow}{q}
\newcommand{\mflow}{\massflow}

\newcommand{\Power}{P}% \power exists in package SIunits
\newcommand{\pressure}{p}
\newcommand{\press}{\pressure}

\newcommand{\temperature}{T}
\newcommand{\temp}{\temperature}
\newcommand{\velocity}{v}
\newcommand{\vel}{\velocity}

\DeclareSIUnit{\kWh}{kWh}
\newcommand{\depot}{\text{d}}
\newcommand{\discrEstimate}{\eta^{\text{d}}}
\newcommand{\modelEstimate}{\eta^{\text{m}}}
\newcommand{\errorest}{\eta}
\newcommand{\errorestd}{\eta^{\text{d}}}
\newcommand{\errorestm}{\eta^{\text{m}}}
\newcommand{\errorex}{\nu}
\newcommand{\errorexd}{\nu^{\text{d}}}
\newcommand{\errorexm}{\nu^{\text{m}}}
\newcommand{\errortol}{\varepsilon}
\newcommand{\errortune}{\tau}

\newcommand{\upswitchingParam}{\varTheta_\mathcal{U}}
\newcommand{\refinementParam}{\varTheta_\mathcal{R}}
\newcommand{\coarseningParam}{\varTheta_\mathcal{C}}
\newcommand{\downswitchingParam}{\varTheta_\mathcal{D}}

\newcommand{\safeguard}{\mu}

\makeatletter
\newenvironment{varsubequations}[1]
 {%
  \addtocounter{equation}{-1}%
  \begin{subequations}
  \renewcommand{\theparentequation}{#1}%
  \def\@currentlabel{#1}%
 }
 {%
  \end{subequations}\ignorespacesafterend
 }
\makeatother

\let\oldref\refeq
\renewcommand{\refeq}[1]{(\oldref{#1})}

\newcommand{\pipes}{\Apipe}
\newcommand{\refine}{\mathcal R}
\newcommand{\coarsen}{\mathcal C}
\newcommand{\upmodel}{\mathcal U}
\newcommand{\downmodel}{\mathcal D}
\newcommand{\dotle}{\mathrel{\dot \le}}

\newcommand{\rev}[1]{#1}

%% file: abstract.tex
We propose an adaptive optimization algorithm for operating district
heating networks in a stationary regime.
The behavior of hot water flow in the pipe network is modeled using
the incompressible Euler equations and a suitably chosen energy
equation.
By applying different simplifications to these equations, we derive a
catalog of models.
Our algorithm is based on this catalog and adaptively controls
where in the network which model is used.
Moreover, the granularity of the applied discretization is controlled
in a similar adaptive manner.
By doing so, we are able to obtain optimal solutions at low
computational costs that satisfy a prescribed tolerance w.r.t.\ the
most accurate modeling level.
To adaptively control the switching between different levels and the
adaptation of the discretization grids, we derive error \rev{measure} formulas
and a posteriori error \rev{measure} estimators.
Under reasonable assumptions we prove that the adaptive algorithm
terminates after finitely many iterations.
Our numerical results show that the algorithm is able to produce
solutions for problem instances that have not been solvable before.

%%% Local Variables:
%%% mode: latex
%%% TeX-master: "dhn-adaptive"
%%% End:

%% file: keywords.tex
District heating networks,
Adaptive methods,
Nonlinear optimization%
%
%
%%% Local Variables:
%%% mode: latex
%%% TeX-master: "dhn-adaptive"
%%% End:

%% file: msc2010.tex
90-XX, % Operations research, mathematical programming
90Cxx, % Mathematical programming
90C11, % Mixed integer programming
90C35, % Programming involving graphs or networks
90C90% Applications of mathematical programming
%
%%% Local Variables:
%%% mode: latex
%%% TeX-master: "dhn-adaptive"
%%% End:

%% file: introduction.tex
\section{Introduction}
\label{sec:introduction}

An efficient and sustainable energy sector is at the core of the
fight against the climate crisis.
Thus, many countries around the world strive towards an energy
turnaround with the overarching goal to replace fossil fuels with
energy from renewable resources such as wind and solar power.
However, one then faces issues with the high volatility of the
fluctuating renewable resources.
To overcome this fluctuating nature of wind and solar power, two main
approaches are currently seen as the most promising ones:
(i) the development and usage of large-scale energy storage systems as
well as (ii) sector-coupling.

In this paper, we consider the computation of optimal operation
strategies for district heating networks.
These networks are used to provide customers with hot water in order
to satisfy their heat demand.
Thus, a district heating network can be seen both as a large-scale
energy storage as well as a key element of successful sector-coupling.
The hot water in the pipes of a district heating network is heated in
so-called depots in which, usually, waste incineration is used as
the primary heat source.
If, however, waste incineration is not sufficient for heating the
water, gas turbines are used as well.
The hot water in the pipeline system can thus be seen as an energy
storage that could, for instance, also be filled using power-to-heat
technologies in time periods with surplus production of renewables.
On the other hand, heat-to-power can be used to smooth the
fluctuating nature of renewables in time periods with only small
renewable production.
Consequently, district heating networks can be seen as sector-coupling
entities with inherent storage capabilities.

To make such operational strategies for district heating networks
possible, an efficient control of the network is required that does
not compromise the heat demand of the households that are connected
to the network.
However, a rigorous physical and technical modeling of hot water flow
in pipes leads to hard mathematical optimization problems.
At the core of these problems are partial differential equations for
modeling both water and heat transport.
Additionally, proper models of the depot and the households further
increase the level of nonlinearity in the overall model.
Finally, the tracking of water temperatures across nodes of the
network leads to nonconvex and nonsmooth mixing models that put a
significant burden on today's state-of-the-art optimization
techniques.

In this paper, we consider the simplified setting of a stationary flow
regime.
For closed-loop control strategies for instationary variants of the
problem we refer to
\cite{Sandou_et_al:2005,Verrilli_et_al:2017,Benonysson_et_al:1995} and
to \cite{Krug_et_al:2019} for open-loop optimization approaches.
Interestingly, the literature on mathematical optimization for
district heating networks is rather sparse.
An applied case study for a specific district heating network in
South Wales is done in \cite{Pirouti_et_al:2013} and
\cite{Rezaie_Rosen:2012} provides more a general discussion of
technological aspects and the potentials of district heating
networks.
In \cite{Schweiger_et_al:2017}, the authors follow a
first-discretize-then-optimize approach for the underlying
PDE-constrained problem.
For the relation between district heating networks and energy storage
aspects we refer to
\cite{Colella_et_al:2012,Verda_Colella:2011,Ben_Hassine_Eicker:2013}
and the references therein.
Stationary models of hot water flow are also considered in studies on
the design and expansion of networks as, \eg, in
\cite{Roland_Schmidt:2020,Bordin_et_al:2016,Dorfner_Hamacher:2014,Roland_Schmidt:2020}.
Numerical simulation of district heating networks using a local time
stepping method is studied in \cite{Borsche_et_al:2018} and model
order reduction techniques for the hyperbolic equations in district
heating networks are discussed in \cite{Rein_Mohring_et_al:2019}
or~\cite{Rein_et_al:2018,Rein_et_al:2019}.
Finally, a port-Hamiltonian modeling approach for district heating
networks is presented and discussed in \cite{Hauschild_et_al:2020}.

Despite the mentioned simplification of considering stationary flow
regimes, the optimization problems at
hand are still large-scale and highly nonlinear
mathematical programs with complementarity constraints (MPCCs) that are
constrained by ordinary differential equations (ODEs).
It turns out that these models are extremely hard to solve for
realistic or even real-world district heating networks if they are
presented to state-of-the-art optimization solvers.
Our contribution is the development of an adaptive optimization
algorithm that controls the modeling and the discretization of the hot
water flow equations in the network.
A similar approach has already been developed and tested for natural
gas networks in \cite{Mehrmann_et_al:2018}.
The main rationale is that simplified (and thus computationally
cheaper) models can lead to satisfactory (w.r.t.\ their physical
accuracy) results for some parts of the network whereas other parts
require a highly accurate modeling to obtain the required physical
accuracy.
The problem, however, is that it is not known up-front where which
kind of modeling is appropriate.
Our adaptive algorithm is based on (i) a catalog of different models of hot
water flow and on (ii) different discretization grids for the
underlying differential equations.
The proposed method then controls the choice of the model and the
discretization grid separately for every pipe in the network.
The switching between different models and discretization grids is based on
rigorous error measures so that we obtain a finite termination proof stating
that the method computes a locally optimal point that is feasible
w.r.t.\ the most accurate modeling level and a prescribed tolerance.
Besides these theoretical contributions, we also show the
effectiveness of our approach in practice and, in particular,
illustrate that instances on realistic networks can be solved with the
newly proposed method that have been unsolvable before.

The remainder of the paper is structured as follows.
In Section~\ref{sec:problem-statement} we present our modeling of
district heating networks and derive the modeling catalog for hot
water flow as well as the discretizations of the respective
differential equations.
After this, we derive exact error \rev{measures} and error \rev{measure} estimators in
Section~\ref{sec:error} both for modeling as well as discretization
error \rev{measures}.
These are then used in Section~\ref{sec:adaptive} to set up the adaptive
optimization algorithm and to prove its finite termination.
The algorithm is numerically tested in
Section~\ref{sec:numerical-results} before we close the paper with some
concluding remarks and some aspects of potential future work in
Section~\ref{sec:conclusion}.

%%% Local Variables:
%%% mode: latex
%%% TeX-master: "dhn-adaptive"
%%% End:

%% file: modeling.tex
\section{Modeling}
\label{sec:problem-statement}

We model the district heating network as a directed and connected
graph $G = (\Vertices, \Arcs)$, which has a special structure.
First, we have a so-called forward-flow part of the network, which is
used to provide the consumers with hot water.
Second, the cooled water is transported back to the depot in the
so-called backward-flow part.
These two parts are connected via the depot in which the cooled water
is heated again, and via the consumers who use the temperature
difference to satisfy the thermal energy demand in the corresponding
household.
The set of nodes of the forward-flow part is denoted by~$\Vff$ and the
set of arcs of this part is denoted by~$\Aff$, \ie, $\arc = (\node,
\otherNode) \in \Aff$ implies $\node, \otherNode \in \Vff$.
In analogy, the set of nodes of the backward-flow part is denoted by~$\Vbf$
and the set of arcs of this part is denoted by~$\Abf$, \ie, $\arc =
(\node, \otherNode) \in \Abf$ implies $\node, \otherNode \in \Vbf$.
The depot arc is denoted by~$\arc_\depot = (\node, \otherNode)$ with
$\node \in \Vbf$ and $\otherNode \in \Vff$.
The consumers are modeled with arcs~$\arc = (\node,
\otherNode)$ with $\node \in \Vff$ and $\otherNode \in \Vbf$.
Finally, all pipes of the forward and the backward flow part are
contained in the set of pipes~$ \pipes = \Aff \cup \Abf $.

In the next subsection we present the model for all components of the
network, \ie, for pipes, consumers, and the depot.

%%%%%%%%%%%%%%%%%%%%%%%%%%%%%%%%%%%%%%%%%%%%%%%%%%%%%%%%%%%%%%%%%%%%%%%%%%%
\subsection{Pipes}
\label{sec:pipe-modeling}

We now derive an approximation for the stationary water flow in
cylindrical pipes.
This derivation is based on the 1-dimensional compressible Euler
equations~\mbox{\cite{Borsche_et_al:2018, Rein_et_al:2018,
    Hauschild_et_al:2020}}
\begin{subequations}
  \label{eq:euler-equations}
  \begin{align}
    0 &= \frac{\partial \rho_\arc}{ \partial t} + \vel_\arc \frac{\partial
        \rho_\arc}{ \partial x} + \rho_\arc \frac{\partial \vel_\arc}{ \partial
        x},
    \label{eq:euler-equations:cont}\\
    0 &= \frac{\partial (\rho_\arc \vel_\arc)}{\partial t} + \vel_\arc
        \frac{\partial (\rho_\arc \vel_\arc)}{\partial x} + \frac{\partial \press_\arc}{\partial
        x} + \frac{\lambda_\arc}{2 \diam_\arc} \rho_\arc \abs{\vel_\arc}
        \vel_\arc + \grav \rho_\arc h_\arc'.
        \label{eq:euler-equations:mom}
  \end{align}
\end{subequations}
Equation~\eqref{eq:euler-equations:cont} is the continuity equation
and models mass balance, whereas the pressure gradient is described by
the momentum equation~\eqref{eq:euler-equations:mom}.
Here and in what follows, $\rho$ denotes the density of water, $\vel$
its velocity, and $\press$ its pressure.
In~\eqref{eq:euler-equations}, the quantities are to be seen as
functions in space ($x$) and time ($t$), \ie, for instance, $\press
= \press(x,t)$.
The diameter of a pipe~$\arc$ is denoted by~$\diam_\arc$,
$\lambda_\arc$ is the pipe's friction coefficient, and $h_\arc'$
denotes the slope of the pipe.
Finally, $g$ is  the gravitational acceleration.

The incompressibility of water is modeled as $ 0 = \rho_{\arc}
\frac{\partial \vel_\arc}{\partial x}$, \cf\
\cite{Hauschild_et_al:2020}, which implies
\begin{equation}
  \label{eq:incompressibility}
  0 = \frac{\partial \rho_\arc}{ \partial t} + \vel_\arc
  \frac{\partial \rho_\arc}{ \partial x}.
\end{equation}
Moreover, the additional PDEs
\begin{subequations}
  \label{eq:internal-energy-and-entropy}
  \begin{align}
    0 &= \frac{\partial \densintene_\arc}{\partial t} + \vel_\arc
        \frac{\partial \densintene_\arc}{\partial x} +\press_\arc
        \frac{\partial \vel_\arc}{\partial x} - \frac{\lambda_\arc}{2
        \diam_\arc} \rho_\arc \abs{\vel_\arc} \vel_\arc^2 + \frac{4
        \heattrans}{\diam_\arc}(\temp_\arc - \soiltemp),
        \label{eq:internal-energy-and-entropy:energy}\\
    0 &= \frac{\partial \densintent_\arc}{\partial t} + \vel_\arc
        \frac{\partial \densintent_\arc}{\partial x} +
        \frac{\lambda_\arc \rho_\arc }{2 \diam_\arc \temp_\arc}
        \abs{\vel_\arc} \vel_\arc^2 + \frac{4 \heattrans}{\diam_\arc}\frac{(\temp_\arc -
        \soiltemp)}{\temp_\arc}
        \label{eq:internal-energy-and-entropy:entropy}
  \end{align}
\end{subequations}
model conservation of internal energy density~$\densintene$ and
entropy density~$\densintent$, respectively; see
\cite{Hauschild_et_al:2020}.
The water's temperature is denoted by $\temp_\arc$.
The parameters $\heattrans$ and $\soiltemp$ are the heat transfer
coefficient and the soil or pipe wall temperature.

Since we expect the change (in time) of the pressure energy and
the term of energy and power loss due to dissipation work to be
small, we neglect these terms. However, if these terms are taken into
account, then it is possible to reformulate these equations in a
port-Hamiltonian form (see, \eg, \cite{Hauschild_et_al:2020}), which
is more appropriate for sector coupling \cite{MehM19}. Finally, we
are interested in the stationary state of the network. This is
modeled by setting all partial derivatives w.r.t.\ time
to zero. Hence, the
System~\eqref{eq:euler-equations}--\eqref{eq:internal-energy-and-entropy}
simplifies to the stationary, incompressible, and 1-dimensional Euler
equations for hot water pipe flow, \ie,
\begin{subequations}
  \label{eq:stationary-ode-system}
  \begin{align}
    0 &=  \dens_\arc \frac{\diff\vel_\arc}{ \diff x},
        \label{eq:stationary-ode-system:incompressibility}\\
    0 &=\vel_\arc \frac{\diff\rho_\arc}{ \diff x} + \rho_\arc
        \frac{\diff\vel_\arc}{ \diff x},
        \label{eq:stationary-ode-system:mass-conservation} \\
    0 &=  \vel_\arc \frac{\diff(\rho_\arc \vel_\arc)}{\diff x} +
        \frac{\diff\press_\arc}{\diff x} + \frac{\lambda_\arc}{2
        \diam_\arc} \rho_\arc \abs{\vel_\arc} \vel_\arc + \grav
        \rho_\arc h_\arc',  \\
    0 &=  \vel_\arc \frac{\diff\densintene_\arc}{\diff x} +\press_\arc
        \frac{\diff\vel_\arc}{\diff x} - \frac{\lambda_\arc}{2
        \diam_\arc} \rho_\arc \abs{\vel_\arc} \vel_\arc^2 + \frac{4
        \heattrans}{\diam_\arc}(\temp_\arc - \soiltemp),  \\
    0 &=  \vel_\arc \frac{\diff\densintent_\arc}{\diff x} +
        \frac{\lambda_\arc \rho_\arc }{2 \diam_\arc \temp_\arc}
        \abs{\vel_\arc} \vel_\arc^2 + \frac{4 \heattrans}{\diam_\arc}\frac{(\temp_\arc -
        \soiltemp)}{\temp_\arc}.
  \end{align}
\end{subequations}
Since $\rho_\arc > 0$ holds,
Equation~\eqref{eq:stationary-ode-system:incompressibility} implies
that $\vel_\arc(x) = \vel_\arc$ is constant for all pipes.
Using this, \eqref{eq:stationary-ode-system:mass-conservation}
implies that the density~$\rho_\arc(x) = \rho_\arc$ is constant as
well.
In addition, we set $\rho_\arc = \rho$ for all arcs~$\arc$ of the
network.
With the mass flow
\begin{equation}
  \label{eq:vel-to-mass-flow-conversion}
  \mflow_\arc = \area_\arc \rho \vel_\arc
\end{equation}
and constant velocities and densities we also have that
$\mflow_\arc(x) = \mflow_\arc$ is constant for all pipes.
In~\eqref{eq:vel-to-mass-flow-conversion}, $\area_\arc$ denotes the
cross-sectional area of pipe~$\arc$. By subsuming the discussed
simplifications we get the system
\begin{subequations}
  \begin{align}
    0 &=  \frac{\diff\press_\arc}{\diff x} + \frac{\lambda_\arc}{2
        \diam_\arc} \rho \abs{\vel_\arc} \vel_\arc + \grav
        \rho h_\arc', \label{eq:momentum-conservation}  \\
    0 &=  \vel_\arc \frac{\diff\densintene_\arc}{\diff x} - \frac{\lambda_\arc}{2
        \diam_\arc} \rho \abs{\vel_\arc} \vel_\arc^2 + \frac{4
        \heattrans}{\diam_\arc}(\temp_\arc - \soiltemp), \\
    0 &=  \vel_\arc \frac{\diff\densintent_\arc}{\diff x} +
        \frac{\lambda_\arc \rho }{2 \diam_\arc \temp_\arc}
        \abs{\vel_\arc} \vel_\arc^2 + \frac{4 \heattrans}{\diam_\arc}\frac{(\temp_\arc -
        \soiltemp)}{\temp_\arc}.
  \end{align}
\end{subequations}
In Equation~\eqref{eq:momentum-conservation}, the pressure gradient
term is the only term that depends on the spatial position~$x$.
Hence, we obtain the stationary momentum and energy equation
\begin{varsubequations}{M1}
  \label{eq:model1}
  \begin{align}
    0 &= \frac{\press_\arc(\length_\arc) -
        \press_\arc(0)}{\length_\arc} + \frac{\lambda_\arc}{2
        \diam_\arc} \rho \abs{\vel_\arc} \vel_\arc + \grav
        \rho h_\arc',\\
    0 &= \vel_\arc \frac{\diff\densintene_\arc}{\diff x} -
        \frac{\lambda_\arc}{2 \diam_\arc} \rho \abs{\vel_\arc}
        \vel_\arc^2 + \frac{4
        \heattrans}{\diam_\arc}(\temp_\arc - \soiltemp). \label{eq:model1-energy}
  \end{align}
\end{varsubequations}

In the following, for our optimization framework, we do not consider
the entropy equation, which can be solved in a post-processing step
once the optimal pressure and internal energy values have been
determined.
%is dropped \HD{Again, shouldn't we provide motivation for %dropping
%an ODE?}since it is supposed unnecessary %in our context.
%\JMR{After discussion with Volker and Hannes we decided %that entropy
%  is not needed at all in the optimization. We can %compute it based on
%the solution of the optimization problem. We just use the %mixing from
%the article we talked about/ use very easy rules for the %depot and consumer.}

The system is closed by the state equations
\begin{subequations}
  \begin{align}
    \rho &= \SI{997}{\kilo \gram \per \cubic \meter} ,\\
    \temp_\arc &= \theta_2(\densintene_\arc^*)^2 + \theta_1\densintene^*_\arc+
                 \theta_0, \label{eq:state-temperature}
  \end{align}
\end{subequations}
in which we set
\begin{gather*}
  \densintene_\arc^* \define \frac{\densintene_\arc}{\densintene_0}, \quad
  \densintene_0 \define \SI[parse-numbers=false]{10^{9}}{\joule \per \cubic \meter}, \\
  \theta_2 = \SI[parse-numbers=false]{59.2453}{\kelvin}, \quad
  \theta_1 = \SI[parse-numbers=false]{220.536}{\kelvin}, \quad
  \theta_0 = \SI[parse-numbers=false]{274.93729}{\kelvin}.
\end{gather*}
Equation~\eqref{eq:state-temperature} is known to be a reasonable
approximation for $\densintene_\arc \in [0.2,0.5]$\,\si{\giga \joule
  \per \cubic \meter}, $\temp_\arc \in [323,403]$\,\si{\kelvin}, and $\press_\arc \in
[5,25]$\,\si{\bar}; see, \eg, \cite{Hauschild_et_al:2020}.

%%%%%%%%%%%%%%%%%%%%%%%%%%%%%%%%%%%%%%%%%%%%%%%%%%%%%%%%%%%%%%%%%%%%%%%%%%%
\subsubsection{Model Catalog}
\label{sec:pipe-modeling:catalog}

For the adaptive optimization method developed in this work we employ
a catalog of models. In this catalog, System~\eqref{eq:model1}
represents the highest or first modeling level, \ie, the most accurate
one.

To derive the second modeling level, we neglect the
(small) term $\lambda_\arc / (2 \diam_\arc) \rho \vel_\arc^2
\abs{\vel_\arc}$ and get
\begin{varsubequations}{M2}
  \label{eq:model2}
  \begin{align}
    0 &=  \frac{\press_\arc(\length_\arc) -
        \press_\arc(0)}{\length_\arc} + \frac{\lambda_\arc}{2
        \diam_\arc} \rho \abs{\vel_\arc} \vel_\arc + \grav
        \rho h_\arc',
    \\
    0 &=  \vel_\arc \frac{\diff\densintene_\arc}{\diff x} + \frac{4
        \heattrans}{\diam_\arc}(\temp_\arc - \soiltemp).
        \label{eq:model2-energy}
  \end{align}
\end{varsubequations}
By further assuming that the first term in~\eqref{eq:model2-energy}
dominates the second one, we can neglect the term $4
\heattrans / \diam_\arc (\temp_\arc - \soiltemp)$ and simplify
System~\eqref{eq:model2} to obtain the third level as
\begin{align}
  \label{eq:model3} \tag{M3}
  \begin{split}
    0 &=  \frac{\press_\arc(\length_\arc) -
      \press_\arc(0)}{\length_\arc} + \frac{\lambda_\arc}{2
      \diam_\arc} \rho \abs{\vel_\arc} \vel_\arc + \grav
    \rho h_\arc',
    \\
    0 &=  \densintene_\arc(\length_\arc) - \densintene_\arc(0).
  \end{split}
\end{align}

Considering model catalogs such as the one just developed is a
standard procedure in order to cope with challenging optimization
models; see, e.g., \cite{DomHLMMT21} in the context of gas networks.
Under sufficient regularity assumptions that allow for Taylor
expansions, a detailed perturbation analysis and the dropping of
higher-order terms would lead to a similar catalog.

%%%%%%%%%%%%%%%%%%%%%%%%%%%%%%%%%%%%%%%%%%%%%%%%%%%%%%%%%%%%%%%%%%%%%%%%%%%
\subsubsection{Exact Solution of the Energy Equation}
\label{sec:exact-solut-energy}

The equations~\eqref{eq:model1-energy} and
\eqref{eq:model2-energy} can be solved analytically.
This is done in the following lemma and will be used later to
compute exact error \rev{measures} in our adaptive algorithm.

\begin{lemma}
  \label{lem:energy-M1--solution}
  The differential equation~\eqref{eq:model1-energy}, \ie,
  \begin{equation*}
    0 = \vel_\arc \frac{\diff\densintene_\arc}{\diff x} -
    \frac{\lambda_\arc}{2 \diam_\arc} \rho \abs{\vel_\arc}
    \vel_\arc^2 + \frac{4
      \heattrans}{\diam_\arc}(\temp_\arc - \soiltemp),
  \end{equation*}
  with initial condition
  \begin{equation*}
    \densintene_\arc(0) = \densintene_\arc^0 > 0
  \end{equation*}
  and state equation~\eqref{eq:state-temperature} has the exact solution
  \begin{equation*}
    \densintene_\arc(x)
    =\frac{\sqrt{\beta^2-4 \alpha
        \gamma}}{2\alpha}\frac{1+\exp\left(\frac{x\sqrt{\beta^2-4
            \alpha \gamma}}{\zeta}\right)
      \left( \frac{2\alpha\densintene_\arc^0+\beta-\sqrt{\beta^2-4
            \alpha \gamma}}{2\alpha\densintene_\arc^0+\beta
          +\sqrt{\beta^2-4 \alpha \gamma}}\right)}{1-
      \exp\left(\frac{x\sqrt{\beta^2-4 \alpha
            \gamma}}{\zeta}\right)\left(\frac{2\alpha\densintene_\arc^0+\beta-
          \sqrt{\beta^2-4 \alpha \gamma}}{2\alpha\densintene_\arc^0+\beta
          +\sqrt{\beta^2-4 \alpha \gamma}}\right)}- \frac{\beta}{2\alpha}
  \end{equation*}
  with
  \begin{equation}
    \begin{gathered}
      \label{eq:lemma-parameter-definition}
      \alpha \define - \frac{4 \heattrans \theta_2 }{\diam_\arc
        (\densintene_0)^2},
      \quad
      \beta \define - \frac{4 \heattrans \theta_1
      }{\diam_\arc\densintene_0},
      \quad
      \zeta \define \vel_\arc,
      \\
      \gamma \define \frac{\lambda_\arc}{2 \diam_\arc} \rho
      \abs{\vel_\arc} \vel_\arc^2 - \frac{4 \heattrans}{\diam_\arc}
      (\theta_0 - \soiltemp),
    \end{gathered}
  \end{equation}
  if $4\alpha \gamma-\beta^2 < 0$ is satisfied.
\end{lemma}
\input{proof-lemma-energy-M1--solution-v2.tex}

Let us further note that the condition $4 \alpha \gamma - \beta^2 < 0$
of the last lemma is satisfied for usual pipe parameters.

\begin{corollary}
  \label{lem:energy-M2--solution}
  The differential equation~\eqref{eq:model2-energy}, \ie,
  \begin{equation*}
    0 = \vel_\arc \frac{\diff\densintene_\arc}{\diff x} + \frac{4
      \heattrans}{\diam_\arc}(\temp_\arc - \soiltemp),
  \end{equation*}
  with initial condition
  \begin{equation*}
    \densintene_\arc(0) = \densintene_\arc^0 > 0
  \end{equation*}
  and state equation~\eqref{eq:state-temperature} has the solution
  \begin{equation*}
    \densintene_\arc(x)
    =\frac{\sqrt{\beta^2-4 \alpha
        \gamma}}{2\alpha}\frac{1+\exp\left(\frac{x\sqrt{\beta^2-4
            \alpha \gamma}}{\zeta}\right)
      \left( \frac{2\alpha\densintene_\arc^0+\beta-\sqrt{\beta^2-4
            \alpha \gamma}}{2\alpha\densintene_\arc^0+\beta
          +\sqrt{\beta^2-4 \alpha \gamma}}\right)}{1-
      \exp\left(\frac{x\sqrt{\beta^2-4 \alpha
            \gamma}}{\zeta}\right)\left(\frac{2\alpha\densintene_\arc^0+\beta-
          \sqrt{\beta^2-4 \alpha \gamma}}{2\alpha\densintene_\arc^0+\beta
          +\sqrt{\beta^2-4 \alpha \gamma}}\right)}- \frac{\beta}{2\alpha}
  \end{equation*}
  with
  \begin{equation*}
    \alpha \define - \frac{4 \heattrans \theta_2 }{\diam_\arc
      (\densintene_0)^2},
    \quad
    \beta \define - \frac{4 \heattrans \theta_1 }{\diam_\arc
      \densintene_0},
    \quad
    \gamma \define  - \frac{4 \heattrans}{\diam_\arc} (\theta_0 -
    \soiltemp),
    \quad \zeta \define \vel_\arc,
  \end{equation*}
  if $4\alpha \gamma-\beta^2 < 0$ is satisfied.
\end{corollary}
The proof is analogous to the one of
Lemma~\ref{lem:energy-M1--solution}.
Figure~\ref{fig:internal-energy} shows the exact solution of
\eqref{eq:model1-energy} for a specific pipe.

The exact solutions derived in the last lemma and corollary could, in
principle, be used as constraints in a nonlinear optimization model.
However, the fractions, square roots, and exponential functions would
lead to a very badly posed problem resulting in an extreme
numerical challenge even for state-of-the-art solvers.

\begin{figure}
  \centering
  \resizebox{0.60\textwidth}{!}{\input{figures/pgf/energy_solution_positive_velocity.pgf}}
  \caption{Solution of Equation~\eqref{eq:model1-energy} for positive velocities and
    the parameters $\heattrans = \SI{0.5}{\watt \per \meter \squared \per \kelvin}$,
    $\lambda_\arc = 0.017$,
    $\diam_\arc = \SI{0.107}{\meter}$,
    $\length_\arc = \SI{1000}{\meter}$, and $\soiltemp =
    \SI{278}{\kelvin}$. The units of the internal energy density
    $\densintene_\arc$ and the velocity $\vel_\arc$ are
    \si{\giga \joule \per \cubic \metre} and \si{\meter \per \second},
    respectively.}
  \label{fig:internal-energy}
\end{figure}

%%%%%%%%%%%%%%%%%%%%%%%%%%%%%%%%%%%%%%%%%%%%%%%%%%%%%%%%%%%%%%%%%%%%%%%%%%%
\subsubsection{Discretization}
\label{sec:pipe-modeling:discretization}

In order to solve the optimization problem, we follow the
first-discretize-then-optimize approach and introduce an equidistant
discretization
\rev{
  \begin{equation*}
    \spacepointsset_\arc
    = \Defset{x_k = k \Delta x_a}{k = 0, \dotsc, n_a}
    \quad \text{with} \quad
    \Delta x_a = L_a / n_a
  \end{equation*}}%
of the spatial domain $[0, L_a]$ using the discretization points~$ x_k
\in \spacepointsset_\arc$ \rev{with} $0 = x_0 < x_1 < \dotsb < x_{n_a}
= L_a$ \rev{and} step size~$\Delta x_a = L_a / n_a = x_{k+1} - x_k$
for $ k = 0,1, \dotsc, n_a - 1$.
We use the implicit mid-point rule to discretize the separate levels of
the catalog, \ie, Systems~\eqref{eq:model1}--\eqref{eq:model3},
as well as the state equation \eqref{eq:state-temperature}. Using the
abbreviation $\densintene_\arc^k \define \densintene_\arc (x_k)$, we
obtain the discretized system
\begin{align}
  \label{eq:discretized1} \tag{D1}
  \begin{split}
    0 &= \frac{\press_\arc(\length_\arc) -
      \press_\arc(0)}{\length_\arc} + \frac{\lambda_\arc \rho}{2
      \diam_\arc}  \abs{\vel_\arc} \vel_\arc + \grav
    \rho h_\arc', \\
    0 &= \vel_\arc \left( \frac{\densintene_\arc^k -
    \densintene_\arc^{k-1}}{\Delta x_a} \right) - \frac{\lambda_\arc \rho}{2
      \diam_\arc} \abs{\vel_\arc} \vel_\arc^2 + \frac
    {4\heattrans}{\diam_\arc} \left(
      \temp_\arc(\densintene_\arc^k,\densintene_\arc^{k-1}) -
      \soiltemp \right)
    \end{split}
\end{align}
for all $k = 1, \dotsc, n_a$. Discretizing \eqref{eq:model2} analogously leads to
\begin{align}
  \label{eq:discretized2} \tag{D2}
  \begin{split}
    0 &= \frac{\press_\arc(\length_\arc) -
      \press_\arc(0)}{\length_\arc} + \frac{\lambda_\arc \rho}{2
      \diam_\arc}  \abs{\vel_\arc} \vel_\arc + \grav
    \rho h_\arc',\\
    0 &= \vel_\arc \left( \frac{\densintene_\arc^k -
    \densintene_\arc^{k-1}}{\Delta x_a} \right)  + \frac
    {4\heattrans}{\diam_\arc} \left(
      \temp_\arc(\densintene_\arc^k,\densintene_\arc^{k-1})
      -\soiltemp \right)
    \end{split}
\end{align}
for all $k = 1, \dotsc, n_a$. The discretized systems
\eqref{eq:discretized1} and \eqref{eq:discretized2} are closed by the
discretized version of the state
equation~\eqref{eq:state-temperature}, \ie, by
\begin{equation}
  \temp_\arc(\densintene_\arc^k,\densintene_\arc^{k-1}) \define
  \frac{\theta_2}{4 \densintene_0^2} \left( \densintene_\arc^k +
    \densintene_\arc^{k-1}\right)^2 +
  \frac{\theta_1}{2\densintene_0} \left( \densintene_\arc^k +
    \densintene_\arc^{k-1}\right) + \theta_0
  \label{eq:state-temperature-discr}
\end{equation}
for all $k = 1, \dotsc, n_a$.
For System \eqref{eq:model3}, we get
\begin{align}
  \label{eq:discretized3} \tag{D3}
  \begin{split}
    0 &= \frac{\press_\arc(\length_\arc) -
      \press_\arc(0)}{\length_\arc} + \frac{\lambda_\arc \rho}{2
      \diam_\arc}  \abs{\vel_\arc} \vel_\arc + \grav
    \rho h_\arc',\\
    0 &= \densintene_\arc^k - \densintene_\arc^{k-1},
  \end{split}
\end{align}
for all $k = 1, \dotsc, n_a$.
In our actual computations, we replace the $n_a$ equations
for~$\densintene$ by the single constraint
$\densintene_\arc(\length_\arc) = \densintene_\arc(0)$, since
a two-point discretization is always exact for this model level.
\begin{figure}
  \centering
  \includegraphics{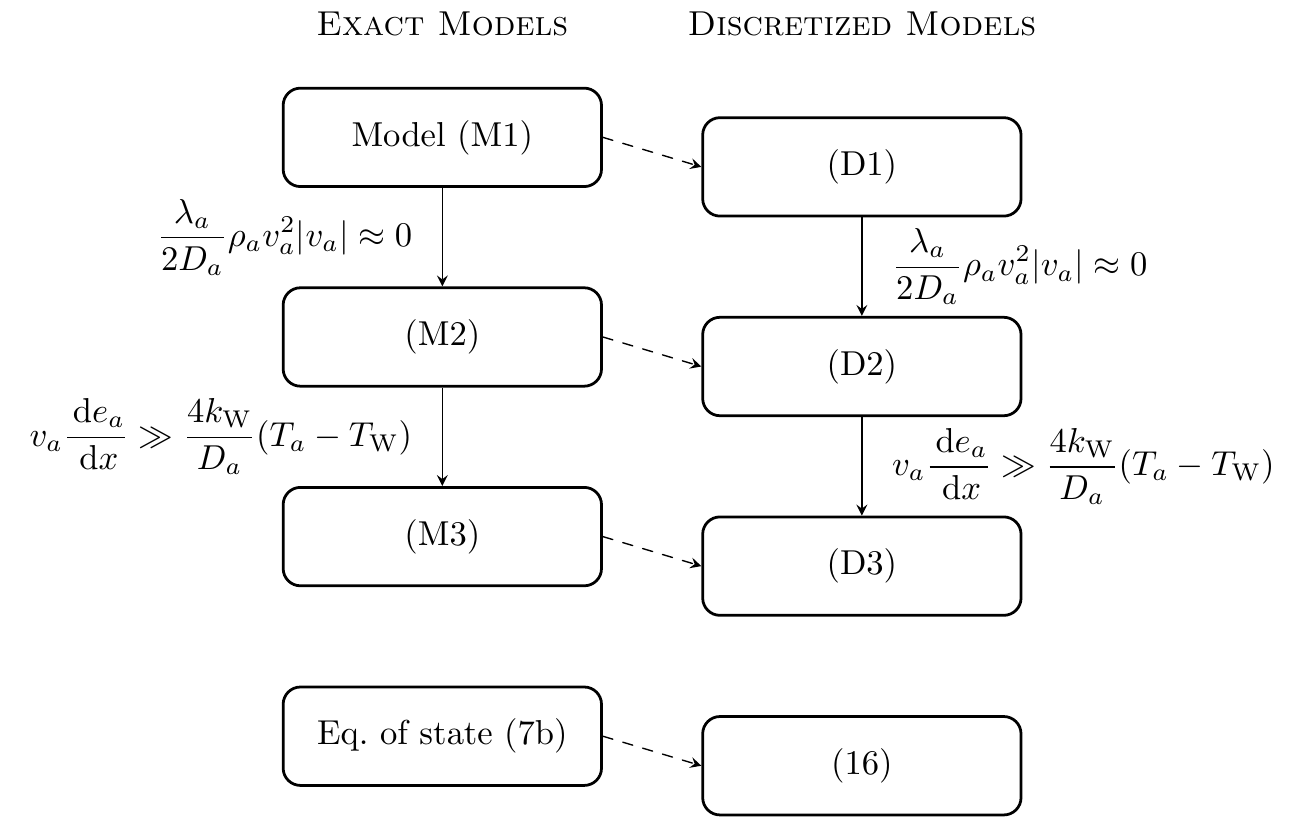}
  \caption{Model catalog for hot water flow in a pipe.}
  \label{fig:summary-adaptive}
\end{figure}
The model catalog (both for the original and the discretized
version) is depicted in Figure~\ref{fig:summary-adaptive}.

%%%%%%%%%%%%%%%%%%%%%%%%%%%%%%%%%%%%%%%%%%%%%%%%%%%%%%%%%%%%%%%%%%%%%%%%%%%
\subsection{Nodes}
\label{sec:nodes}

In this section, we discuss the modeling of nodes in the district
heating network.
To this end, we mainly follow the modeling approach presented
in~\cite{Krug_et_al:2019}.
We model mass conservation via
\begin{equation}
  \label{eq:mass-flow}
  \sum_{\arc \in \Inedges(\node)} \mflow_{\arc}
  = \sum_{\arc \in \Outedges(\node)} \mflow_{\arc},
  \quad \node \in \Vdh,
\end{equation}
where $\Inedges(\node)$ and $\Outedges(\node)$ model the set of in-
and outgoing arcs of node~$\node$, respectively.
We assume continuity of pressure at the nodes and obtain
\begin{align}
  \label{eq:pressure-continuity}
  \begin{split}
    \press_\node
    &= \press_\arc (0), \quad \node \in \Nodes, \ \arc \in
    \Outedges(\node),
    \\
    \press_\node
    &= \press_\arc (\length_\arc), \quad \node \in \Nodes, \ \arc \in
    \Inedges(\node),
  \end{split}
\end{align}
where~$\press_\node$ is the pressure at node~$\node$.

Finally, we have to model how the internal energy is mixed at the
nodes of the network.
%In contrast to the continuity of pressures\HD{say that %this is an
%assumption?}, internal energy density is
To describe this, we use the perfect mixing model
\begin{subequations}
   \label{eq:temperature-mixing}
  \begin{align}
    \label{eq:distr-heat-temperature-mixing-in-complementary}
     \sum_{\arc \in \Inedges(\node)}
    \frac{\densintene_{\arc:\node}\mflow_{\arc}}{\rho}
    &= \sum_{\arc \in \Outedges(\node)}
      \frac{\densintene_{\arc:\node}\mflow_{\arc}}{\rho} ,
    \\
    \label{eq:distr-heat-temperature-mixing-out-complementary-out}
    0 & = \posmflow_{\arc}(\densintene_{\arc:\node} - \densintene_\node) ,
    \quad \arc \in \Outedges(\node),
    \\
    \label{eq:distr-heat-temperature-mixing-out-complementary-in}
    0 & = \negmflow_{\arc}(\densintene_{\arc:\node} - \densintene_\node) ,
    \quad \arc \in \Inedges(\node),
  \end{align}
\end{subequations}
with
\begin{equation}
  \label{eq:mpcc-max-reform}
  \mflow_{\arc} = \posmflow_{\arc} -
  \negmflow_{\arc},
  \quad \posmflow_{\arc} \geq 0,
  \quad \negmflow_{\arc} \geq 0,
  \quad \posmflow_{\arc} \negmflow_{\arc} = 0
\end{equation}
for $\arc \in \nodeArcs(\node) = \Inedges(\node) \cup \Outedges(\node)$.
Here and in what follows, we denote with $e_{a:u}$ the internal energy
in pipe~$a$ at its end node~$u$.
For more details and a derivation of this model we refer to
\cite{Schmidt_et_al:2016,Krug_et_al:2019,Hauschild_et_al:2020,Hante_Schmidt:2019a}.

%%%%%%%%%%%%%%%%%%%%%%%%%%%%%%%%%%%%%%%%%%%%%%%%%%%%%%%%%%%%%%%%%%%%%%%%%%%
\subsection{Depot and Consumers}
\label{sec:depot-consumers}

Following \cite{Krug_et_al:2019,Roland_Schmidt:2020}, for the
depot~$\arc_\depot = (\node, \otherNode)$ we have the constraints
\begin{subequations}%
  \label{eq:depot}%
  \begin{align}
    \press_{\node}
    & = \press_\stagnation,
    \\
    \Powerpress
    & =
      \frac{\mflow_{\arc_\depot}}{\dens} \left(
      \press_{\arc_\depot:\otherNode} -
      \press_{\arc_\depot:\node}\right),
    \\
    \Powerwaste + \Powergas
    & =    \frac{\mflow_{\arc_\depot}}{\dens} \left(
      \densintene_{\arc_\depot:\otherNode} -
      \densintene_{\arc_\depot:\node}\right),
  \end{align}
\end{subequations}
where $\press_\stagnation$ is the so-called stagnation pressure that
is used to make the overall pressure profile in the network unique.
Moreover, $\Powerpress$ is the power required for the pressure
increase realized at the depot, $\Powerwaste$ is the power obtained by
waste incineration and~$\Powergas$ is the power obtained by burning natural gas.
The latter two quantities are used in the depot to increase the internal energy
density or, equivalently, the temperature of the water.

In order to model the consumers~$\arc = (\node, \otherNode) \in
\Acons$, we use the constraints
\begin{subequations}
  \label{eq:consumer}
  \begin{align}
    \Power_\arc
    & = \frac{\mflow_{\arc}}{\dens} \left( \densintene_{\arc:\otherNode}
      - \densintene_{\arc:\node}\right),
      \label{eq:power-demand}\\
    \densintene_{\arc:\node}
    & \geq \densinteneffarc,
      \label{eq:temperature-foreflow}\\
    \densintene_{\arc:\otherNode}
    & = \densintenebf,
      \label{eq:temperature-backflow}\\
    \pressure_{\otherNode}
    & \leq \pressure_{\node}.
    \label{eq:pressure-difference}
  \end{align}
\end{subequations}
The first constraint models how the required thermal
energy~$\Power_\arc$ is obtained in dependence on the  mass
flow~$\mflow_{\arc}$ at the consumer and the
difference~$\densintene_{\arc:\otherNode} - \densintene_{\arc:\node}$
of the internal energy density.
The internal energy density at inflow conditions
($\densintene_{\arc:\node}$) needs to be larger than the given
threshold~$\densinteneffarc$ and, at outflow conditions, it is fixed to
the network-wide constant~$\densintenebf$.
Finally, the fourth constraint states that the pressure cannot be
increased at the household of a consumer.

%%%%%%%%%%%%%%%%%%%%%%%%%%%%%%%%%%%%%%%%%%%%%%%%%%%%%%%%%%%%%%%%%%%%%%%%%%%
\subsection{Bounds, Objective Function, and Model Summary}
\label{sec:bounds}

To complete the model we need to incorporate some technical and
physical bounds on the variables of the model and to define a proper
objective function.
First, we have bounds on the mass flow, \ie,
\begin{equation}
  \label{eq:mass-flow-bounds}
  \mflow_\arc^- \leq \mflow_\arc \leq \mflow_\arc^+,
  \quad \arc \in \Aff \cup \Abf \cup \Acons,
  % \label{eq:mass-flow-bounds-candidate}
\end{equation}
on the nodal pressures,
\begin{equation}
  \label{eq:pressure-bounds}
  0  \leq \press_\node
  \leq \press_\node^+,
  \quad
  \node \in \Vdh,
\end{equation}
and on the nodal water temperatures, \ie,
\begin{equation}
  \label{eq:temperature-bounds}
  \temp_\node \in [\temp_\node^-,\temp_\node^+],
  \quad \node \in \Vdh.
\end{equation}
Lastly, we incorporate bounds on power consumption, \ie,
\begin{equation}
  \label{eq:power-bounds}
  \Powerpress  \in [0,\Powerpress^+],
  \quad
  \Powerwaste  \in [0,\Powerwaste^+],
  \quad
  \Powergas  \in [0,\Powergas^+]
\end{equation}
for given upper bounds~$\Powerpress^+$, $\Powerwaste^+$, and
$\Powergas^+$.

Our goal is to minimize the overall costs
required to satisfy the heat demand of all the consumers.
Thus, the objective function is given by
\begin{equation}
  \label{eq:objective}
  \pressurecost\Powerpress + \wastecost\Powerwaste + \gascost\Powergas,
\end{equation}
where $\pressurecost,\wastecost$, and $\gascost$, respectively,
correspond to the cost of pressure increase, waste incineration, and
burning gas.

Taking this all together leads to the discretized and thus
finite-dimensional optimization problem
\begin{equation}
  \label{eq:overall-problem}
  \tag{NLP}
  \begin{split}
    \min \quad & \text{objective: } \eqref{eq:objective},\\
    \st \quad & \text{pipe flow and thermal modeling:
      \eqref{eq:state-temperature-discr} and \eqref{eq:discretized1},
      \eqref{eq:discretized2}, or \eqref{eq:discretized3},}\\
    & \text{mass conservation: \eqref{eq:mass-flow},}\\
    & \text{pressure continuity: \eqref{eq:pressure-continuity},}\\
    & \text{temperature mixing: \eqref{eq:temperature-mixing},
      \eqref{eq:mpcc-max-reform},}\\
    & \text{depot constraints: \eqref{eq:depot},}\\
    & \text{consumer constraints: \eqref{eq:consumer},}\\
    & \text{bounds: \eqref{eq:mass-flow-bounds}--\eqref{eq:power-bounds}.}
  \end{split}
\end{equation}
This is a highly nonlinear and, depending on the chosen grids,
large-scale optimization problem.
Moreover, it only possesses very few degrees of freedom since almost
all variables are determined by our physical modeling.
Both aspects already make the problem very challenging to solve.
In addition, however, the model also contains the complementarity
constraints~\eqref{eq:mpcc-max-reform}, which makes it an
ODE-constrained mathematical program with complementarity constraints (MPCC).
Solving it for real-world networks is very challenging, which is the
motivation of the error \rev{measure}-based adaptive algorithm that we develop in
the two following sections.

%%% Local Variables:
%%% mode: latex
%%% TeX-master: "dhn-adaptive"
%%% End:

%% file: proof-lemma-energy-M1--solution-v2.tex
\begin{proof}
  We combine \eqref{eq:model1-energy} and \eqref{eq:state-temperature}
  to obtain
  \begin{equation*}
    0 =  \vel_\arc \frac{\diff\densintene_\arc}{\diff x} - \frac{\lambda_\arc}{2
      \diam_\arc} \rho \abs{\vel_\arc} \vel_\arc^2 + \frac{4
      \heattrans}{\diam_\arc}(\theta_2(\densintene_\arc^*)^2 +
    \theta_1\densintene^*_\arc + \theta_0 - \soiltemp).
  \end{equation*}
  After re-organizing and replacing $\densintene_\arc^*$ by its
  definition, the equation reads
  \begin{equation}
    - \frac{4 \heattrans \theta_2}{\diam_\arc (\densintene_0)^2}(\densintene_\arc)^2 -
    \frac{4 \heattrans \theta_1}{\diam_\arc \densintene_0}
    \densintene_\arc
    - \frac{4 \heattrans}{\diam_\arc} (\theta_0 - \soiltemp)
    + \frac{\lambda_\arc}{2 \diam_\arc} \rho \abs{\vel_\arc}
    \vel_\arc^2
    =  \vel_\arc \frac{\diff\densintene_\arc}{\diff
      x}. \label{energy-with-state}
  \end{equation}
  We combine Equation~\eqref{energy-with-state} with the definitions
  in~\eqref{eq:lemma-parameter-definition} and get
  \begin{equation}
    \alpha \densintene_\arc^2 + \beta\densintene_\arc  + \gamma =  \zeta
    \frac{\diff\densintene_\arc}{\diff x}. \label{riccati-ode}
  \end{equation}
  Equation~\eqref{riccati-ode} is a \rev{special type of} Riccati
  equation with constant coefficients; see, \eg, \cite{Reid:1972}.
  Because $\alpha,\beta,\gamma$, and $\zeta$ do
  not depend on $x$, they can be seen as constants when integrating over
  $x$. We re-organize and integrate both sides over $x$, yielding
  \begin{equation}
    \int 1 \diff x=\int \frac{\zeta\frac{\diff\densintene_\arc}{\diff
        x}}{\alpha(\densintene_\arc)^2 + \beta\densintene_\arc  +
      \gamma} \diff x. \label{new-riccati-ode}
  \end{equation}
  Applying a variable change in the right-hand side of
  \eqref{new-riccati-ode} leads to
  \begin{equation}
    \frac{x}{\zeta}=\int \frac{1}{\alpha(\densintene_\arc)^2 + \beta\densintene_\arc  +
      \gamma} \diff\densintene_\arc. \label{integral-start}
  \end{equation}
  We may rewrite
  \begin{align*}
    \alpha(\densintene_\arc)^2 + \beta\densintene_\arc  +
    \gamma  &= \left(   (\densintene_\arc)^2 + \frac{\beta}{\alpha}\densintene_\arc +
              \frac{\beta^2}{4\alpha^2}\right) + \frac{4\alpha \gamma - \beta^2}{4\alpha}\\
            &= \alpha \left( \left( \densintene_\arc + \frac{\beta}{2\alpha} \right)^2 +
              \frac{4\alpha \gamma - \beta^2}{4\alpha^2} \right),
  \end{align*}
  since $4\alpha \gamma-\beta^2 < 0$ holds by assumption.
  Therefore, we have
  \begin{equation*}
    \alpha \left( \left( \densintene_\arc + \frac{\beta}{2\alpha} \right)^2 +
      \frac{4\alpha \gamma -
        \beta^2}{4\alpha^2} \right)  = \alpha \left( \left( \densintene_\arc + \frac{\beta}{2\alpha} \right)^2 -
      \left(\frac{\sqrt{\beta^2- 4\alpha \gamma}}{2\alpha}\right)^2 \right).
  \end{equation*}
  Going back to \eqref{integral-start} we have (see also
    Section~8.1 in~\cite{MR3443347})
  \begin{align*}
    &\hspace*{12pt} \int \frac{1}{\alpha(\densintene_\arc)^2 + \beta\densintene_\arc  + \gamma} \diff
      \densintene_\arc \\
    &= \int \frac{1}{\alpha \left( \left(\densintene_\arc +
      \frac{\beta}{2\alpha} \right)^2-
      \left(\frac{\sqrt{\beta^2-4\alpha \gamma}}{2\alpha}\right)^2 \right)}
      \diff\densintene_\arc,\\
    &= C_1 \int \frac{
      \frac{\sqrt{\beta^2- 4\alpha \gamma}}{\alpha}  }{
      \left( \densintene_\arc + \frac{\beta}{2\alpha}
      \right)^2 - \left(\frac{\sqrt{\beta^2-
      4 \alpha \gamma}}{2\alpha}\right)^2 }
      \diff\densintene_\arc, \\
    &= C_1 \int
      \frac{\frac{\sqrt{\beta^2- 4 \alpha \gamma}}{2\alpha} + \frac{\sqrt{\beta^2- 4 \alpha \gamma}}{2\alpha} +\densintene_\arc + \frac{\beta}{2\alpha} -
      \densintene_\arc - \frac{\beta}{2\alpha}}{
      \left( \densintene_\arc + \frac{\beta}{2\alpha}
      \right)^2 - \left(\frac{\sqrt{\beta^2-
      4 \alpha \gamma}}{2\alpha}\right)^2 }
      \diff\densintene_\arc,\\
    &= C_1 \int \left(
      \left(
      \densintene_\arc +
      \frac{\beta}{2\alpha}-\frac{\sqrt{\beta^2-
      4 \alpha \gamma}}{2\alpha} \right)^{-1} -
      \left(
      \densintene_\arc +
      \frac{\beta}{2\alpha}+\frac{\sqrt{\beta^2-
      4 \alpha \gamma}}{2\alpha} \right)^{-1} \right)
      \diff\densintene_\arc ,\\
    &= C_1 \ln \Abs{ \frac{2\alpha \densintene_\arc +
      \beta -\sqrt{\beta^2-
      4 \alpha \gamma}}{2\alpha\densintene_\arc +
      \beta +\sqrt{\beta^2-
      4 \alpha \gamma}} } + C_2,
  \end{align*}
  where we set
  \begin{equation*}
    C_1 \define \frac{1}{\sqrt{\beta^2- 4\alpha \gamma}}.
  \end{equation*}
  The internal energy equation thus reduces to
  \begin{equation*}
    \frac{x}{\zeta}= C_1  \ln \Abs{ \frac{2\alpha
        \densintene_\arc(x) + \beta -\sqrt{\beta^2- 4 \alpha \gamma}}{2a\densintene_\arc(x) + \beta
        +\sqrt{\beta^2-4 \alpha \gamma}} } + C_2.
  \end{equation*}
  By re-substituting the definition of $C_1$ we may write
  \begin{equation}
    \sqrt{\beta^2- 4 \alpha \gamma} \left( \frac{x}{\zeta} - C_2
    \right)
    =   \ln \Abs{ \frac{2\alpha
        \densintene_\arc(x) + \beta -\sqrt{\beta^2- 4 \alpha \gamma}}{2a\densintene_\arc(x) + \beta
        +\sqrt{\beta^2-4 \alpha \gamma}} }. \label{eq:constant-reduction}
  \end{equation}
  We define $C_3 \define \exp \left(-C_2\sqrt{\beta^2- 4 \alpha
      \gamma}\right)$.
  Then, \eqref{eq:constant-reduction} leads to
  \begin{equation*}
    C_3\exp \left( \frac{x \sqrt{\beta^2- 4 \alpha
          \gamma}}{\zeta}\right)= \Abs{ \frac{2\alpha
        \densintene_\arc(x) + \beta -\sqrt{\beta^2- 4 \alpha
          \gamma}}{2a\densintene_\arc(x) + \beta + \sqrt{\beta^2-4
          \alpha \gamma}} }.
  \end{equation*}
  The constant~$C_3$ then absorbs the $\pm$ sign such that we can write
  \begin{equation}
    C_3\exp \left( \frac{x \sqrt{\beta^2- 4 \alpha
          \gamma}}{\zeta}\right)= \left( \frac{2\alpha
        \densintene_\arc(x) + \beta -\sqrt{\beta^2- 4 \alpha
          \gamma}}{2a\densintene_\arc(x) + \beta + \sqrt{\beta^2-4
          \alpha \gamma}} \right). \label{eq:constant-extraction}
  \end{equation}
  We compute $C_2$ using the initial condition at $x=0$ and obtain
  \begin{equation}
    C_3 = \left( \frac{2\alpha
        \densintene_\arc^0 + \beta -\sqrt{\beta^2- 4 \alpha
          \gamma}}{2a\densintene_\arc^0 + \beta + \sqrt{\beta^2-4
          \alpha \gamma}} \right). \label{eq:constant-value}
  \end{equation}
  Finally, we combine Equation \eqref{eq:constant-extraction} and
  \eqref{eq:constant-value}, yielding
  \begin{equation*}
    \densintene_\arc(x)
    =\frac{\sqrt{\beta^2-4 \alpha \gamma}}{2\alpha}\frac{1+\exp\left(\frac{x\sqrt{\beta^2-4 \alpha \gamma}}{\zeta}\right)
      \left( \frac{2\alpha\densintene_\arc^0+\beta-\sqrt{\beta^2-4 \alpha \gamma}}{2\alpha\densintene_\arc^0+\beta
          +\sqrt{\beta^2-4 \alpha \gamma}}\right)}{1-
      \exp\left(\frac{x\sqrt{\beta^2-4 \alpha \gamma}}{\zeta}\right)\left(\frac{2\alpha\densintene_\arc^0+\beta-
          \sqrt{\beta^2-4 \alpha \gamma}}{2\alpha\densintene_\arc^0+\beta
          +\sqrt{\beta^2-4 \alpha \gamma}}\right)}-
    \frac{\beta}{2\alpha}.
    \ifPreprint
    \qedhere
    \fi
  \end{equation*}
\end{proof}

%%% Local Variables:
%%% mode: latex
%%% TeX-master: "dhn-adaptive"
%%% End:

%% file: error-measures.tex
\section{Error Measures}
\label{sec:error}

In this section, we introduce the error measures for the adaptive
optimization algorithm that is presented in
Section~\ref{sec:adaptive}. Our approach is based on the work of
\cite{Mehrmann_et_al:2018} and adapted for the problem at hand. The
algorithm developed here is designed to iteratively solve the
nonlinear program \eqref{eq:overall-problem} until its solution $
\solution $ is deemed to be feasible w.r.t.\ a prescribed tolerance.
The algorithm iteratively switches the model level and the step sizes
of the discretization grids for each pipe according to a switching
strategy presented later on. Both the switching strategy and the
feasibility check utilize the error measures in this section.

For the \eqref{eq:overall-problem}, four error sources can be
identified: errors as introduced by the solver of the optimization
problem, round-off errors, errors from switching between
Systems~\eqref{eq:discretized1}--\eqref{eq:discretized3}, and errors
due to selecting different step sizes of the discretization of the
systems.
In this work we will only consider the latter two error sources, which
we refer to as model (level) error \rev{measures} and discretization
(level) error \rev{measures}, respectively. For a discussion of the neglected
solver and round-off errors we refer to
Remark~\ref{remark:other-error-sources} below.
By investigating the
Systems~\eqref{eq:discretized1}--\eqref{eq:discretized3} one finds
that the only difference between them, and hence the resulting error
source, is the energy equation and its discretization. This is why we
base the definition of the error in each pipe $ a $ on its internal
energy density $ \densintene_\arc $.

In general, utilizing estimates of the error of a system allows for
the assessment of the quality of their solution if an exact solution
is not available. Hence, this section is used to introduce error \rev{measure}
estimates for the model and discretization error \rev{measure}. However, since we
have the analytic solution of the energy equations of
Systems~\eqref{eq:model1}--\eqref{eq:model3} at hand, we can compute
exact error \rev{measures} for the model and discretization error \rev{measure}. Having the exact
error \rev{measures} available allows us to compare them to the error \rev{measure} estimates
presented in this work and, hence, determine their quality.

This section is structured as follows. We start by providing the
required notation in Section~\ref{sec:error:notation}. Furthermore, the
rules that are used to refine and coarsen the grids in the
discretization of
Systems~\eqref{eq:discretized1}--\eqref{eq:discretized3} are
introduced. In Section~\ref{sec:error:measures}, we continue by deriving
exact and estimated error measures. We then close this section by proving
that the error \rev{measure} estimates form upper bounds of the exact error \rev{measures} in a
first-order approximation.

\begin{remark}
  \label{remark:other-error-sources}
  Since we want to be able to employ different third-party optimization
  software packages in our adaptive error control we do not incorporate
  the errors introduced by the solvers for the optimization
  problem. However, if error \rev{measure} estimates and error \rev{measure} control for these
  error \rev{measures} are available then these can  be incorporated as well. It has
  been observed in the application of adaptive methods for gas networks
  \cite{Mehrmann_et_al:2018,Stolwijk_Mehrmann:2018} that round-off
  errors typically do not contribute much to the global error. For this
  reason, we also do not consider round-off errors in our adaptive
  procedure for district heating networks.
\end{remark}

%%%%%%%%%%%%%%%%%%%%%%%%%%%%%%%%%%%%%%%%%%%%%%%%%%%%%%%%%%%%%%%%%%%%%%%%%%%
\subsection{Notation}
\label{sec:error:notation}

We start this section by introducing the required quantities and
notation.
In order to keep the notation lucid, we omit the usage of the subscript~$
\arc $ as much as possible in this section. In particular, we drop the
subscript $ \arc $ for the model level ($ \level_\arc \to \level $), the
grid size ($ \spacediff_\arc \to \spacediff $), and for the set of
gridpoints ($ \spacepointsset_\arc \to \spacepointsset $), if not stated
otherwise.

Let $ \solution $ denote the solution of the optimization problem
\eqref{eq:overall-problem}. For all pipes $ \arc \in \pipes $ it
contains the approximate solution $
\densintene_\arc^\level(x_k; \spacediff) $ for model level $ \level $
(of pipe model (D$\level$)) and step size $ \spacediff $ (of
discretization grid $ \spacepointsset $) at every grid point $ x_k \in
\spacepointsset $, $ k = 1,\ldots, n $. In addition, for a given pipe
$ \arc $ we denote the exact solution of model~(M$\level $),
evaluated at $ x_k \in \spacepointsset $ as $
\densintene_\arc^\level(x_k) $. Furthermore, for the approximate and
exact solutions we also utilize the notion of $
\densintene_\arc^\level(\spacepointsset; \spacediff) \define
(\densintene_\arc^\level(x_1; \spacediff), \ldots,
\densintene_\arc^\level(x_n; \spacediff))^\top $ and $
\densintene_\arc^\level(\spacepointsset) \define
(\densintene_\arc^\level(x_1), \ldots,
\densintene_\arc^\level(x_n))^\top $, respectively.

\if0 %%%%%%%%%%%%%%%%%%%%%%%%%%%%%%%%%%%%%%%%%%%%%%%%%%
\begin{figure}
  \centering
  \begin{tikzpicture}

    \draw[-stealth] (-.5,0) -- (8.5,0) node[anchor=west] {$ \log_2 \Delta x $};
    \draw (0cm,-1mm) node[anchor=north] {$\Delta x_{i+2}$} --++ (0,2mm);
    \draw (2cm,-1mm) node[anchor=north] {$\Delta x_{i+1}$} --++ (0,2mm);
    \draw (4cm,-1mm) node[anchor=north] {$\Delta x_{i}$} --++ (0,2mm);
    \draw (6cm,-1mm) node[anchor=north] {$\Delta x_{i-1}$}--++ (0,2mm);
    \draw (8cm,-1mm) node[anchor=north] {$\Delta x_{i-2}$} --++ (0,2mm);

    \node[] at (20mm,7mm) {$\leftarrow$ refinement: $ \Delta x_i/2 $};
    \node[] at (60mm,7mm) {coarsening: $ 2 \Delta x_i $ $ \rightarrow $};

    % \draw[<-] (-2.5mm,5mm) -- (37.5mm,5mm) node[midway, anchor=south] {grid refinement: $ \Delta x_i/2 $};
    % \draw[->] (42.5mm,5mm) -- (82.5mm,5mm) node[midway, anchor=south] {grid coarsening: $ 2 \cdot \Delta x_i $};

  \end{tikzpicture}
  \caption{Step-size refinement strategy by halving or doubling the
    current step size $\Delta x_i$}
  \label{fig:grid-size-refinement}
\end{figure}
\fi %%%%%%%%%%%%%%%%%%%%%%%%%%%%%%%%%%%%%%%%%%%%%%%%%%%

We continue by defining the grid refinement and coarsening rules. For
a given pipe~$ \arc $, consider a sequence of grids $ \{
\spacepointsset_i \} $, $ i = 0,1,2,\ldots $, with $ \spacepointsset_i
\define \{x_{k_i}\}_{k_i=\rev{0}}^{n_i} $ and $ \spacediff_{i} =
x_{k_{i+1}} - x_{k_i} $ for $ k_i = \rev{0},\ldots, n_i
$. Moreover, we refer to $ \spacepointsset_0 $ as the reference or
evaluation grid. It is defined by a given number of grid points~$ n_0
$ and the corresponding step size $ \Delta x_{0} \define L_a/(n_0 - 1)
$. Given an arbitrary grid $ \spacepointsset_i $, $ i = 0,1,2, \ldots
$, we perform a grid refinement step by halving its step size $
\spacediff_i $ to get $ \spacediff_{i+1} = \spacediff_i/2 $ of the
refined grid $ \spacepointsset_{i+1} $. Conversely, we perform a
grid coarsening step by doubling $ \spacediff_{i} $ of grid $
\spacepointsset_{i} $ to obtain the coarsened grid $
\spacepointsset_{i-1} $ with step size $ \spacediff_{i-1} = 2
\spacediff_{i} $.
%%%Figure~\ref{fig:grid-size-refinement} depicts a visualization of the
%%%grid refinement and coarsening rules.
Performing grid refinement and coarsening this way ensures that for every $ i =
1, 2, \ldots $ it holds that $ \spacepointsset_0 \subset
\spacepointsset_i $. Therefore, providing a fixed number of
grid points~$ n_0 $ enables us to use the reference grid $
\spacepointsset_0 $ as a common evaluation grid for every refinement
and coarsening step.

%%%%%%%%%%%%%%%%%%%%%%%%%%%%%%%%%%%%%%%%%%%%%%%%%%%%%%%%%%%%%%%%%%%%%%%%%%%
\subsection{Derivation of Error Measures}
\label{sec:error:measures}

In the following, we introduce two error measures: exact error \rev{measures} and
error \rev{measure} estimates. To this end, we consider a single pipe $ a \in
\pipes$.
We start by defining the total exact error \rev{measure} as
\begin{equation}
  \label{eq:exact-total-error}
  \errorex_\arc(\solution) \define \|
  \densintene_\arc^1(\spacepointsset_0) -
  \densintene_\arc^\level(\spacepointsset_0; \spacediff_i)
  \|_\infty,
\end{equation}
where we compare the approximate solution of Model (D$\level$) with
grid size $ \spacediff_i $ to the exact solution of
Model~\eqref{eq:model1}.
Note that $\densintene_\arc^\level(\spacepointsset_0;
\spacediff_i)$ is part of the considered solution~$\solution$ and
that the exact error \rev{measure}~$\densintene_\arc^1(\spacepointsset_0)$ can be,
\eg, computed by using the exact formulas given in
Section~\ref{sec:exact-solut-energy}.
Second, we introduce the exact model error \rev{measure} via
\begin{equation}
  \errorexm_\arc(\solution) \define \|
  \densintene_\arc^1(\spacepointsset_0) -
  \densintene_\arc^\level(\spacepointsset_0) \|_\infty\,,
  \label{eq:exact-model-error}
\end{equation}
where we compare the solutions of models (M$\level$)
and~\eqref{eq:model1}.
Next, we define the exact discretization error \rev{measure} as
\begin{equation}
  \label{eq:exact-discretization-error}
  \errorexd_\arc(\solution) \define \|
  \densintene_\arc^\level(\spacepointsset_0) -
  \densintene_\arc^\level(\spacepointsset_0; \spacediff_i) \|_\infty,
\end{equation}
for which we compare the solution of Model~(D$\level$) with grid size $
\spacediff_i $ to the exact solution of Model~(M$\level$). We continue
by introducing error \rev{measure} estimates. The (total) error \rev{measure} estimate is
defined as the sum of a model error \rev{measure} estimate and a
discretization error \rev{measure} estimate. That is,
\begin{equation}
  \label{eq:error-estimate}
  \errorest_\arc(\solution) \define \errorestm_\arc(\solution) + \errorestd_\arc(\solution)
\end{equation}
with the model error \rev{measure} estimate
\begin{equation}
  \label{eq:model-error-estimate}
  \errorestm_\arc(\solution) \define \|
  \densintene_\arc^1(\spacepointsset_0; \spacediff_i) -
  \densintene_\arc^\level(\spacepointsset_0; \spacediff_i) \|_\infty
\end{equation}
and the discretization error \rev{measure} estimate
\begin{equation}
  \label{eq:discretization-error-estimate}
  \errorestd_\arc(\solution) \define \|
  \densintene_\arc^\level(\spacepointsset_0; \spacediff_i) -
  \densintene_\arc^\level(\spacepointsset_0; \spacediff_{i-1})
  \|_\infty.
\end{equation}
The model error \rev{measure} estimate compares two solutions with the same
discretization scheme but different pipe
models~\eqref{eq:discretized1} and~(D$\level$).
On the other hand, the discretization error \rev{measure} estimate compares two
solutions of the same model but with different discretization schemes
as given by the step sizes~$ \spacediff_i $ and~$ \spacediff_{i-1} $.

By considering the definitions
\eqref{eq:exact-total-error}--\eqref{eq:discretization-error-estimate}
one finds the relation
\begin{align}
  \errorex_\arc(\solution)
  &= \| \densintene_\arc^1(\spacepointsset_0) -
    \densintene_\arc^\level(\spacepointsset_0; \spacediff_i) +
    \densintene_\arc^\level(\spacepointsset_0) -
    \densintene_\arc^\level(\spacepointsset_0) \|_\infty \nonumber
  \\
  &\leq \| \densintene_\arc^1(\spacepointsset_0) -
    \densintene_\arc^\level(\spacepointsset_0) \|_\infty\, + \|
    \densintene_\arc^\level(\spacepointsset_0) -
    \densintene_\arc^\level(\spacepointsset_0; \spacediff_i)
    \|_\infty \nonumber\\
  &= \errorexm_\arc(\solution) + \errorexd_\arc(\solution) \dotle
    \errorestm_\arc(\solution) + \errorestd_\arc(\solution) =
    \errorest_\arc(\solution)
    \label{eq:dotle-for-estimator}
\end{align}
for $ \spacediff_{i} \to 0 $.
\input{proof-error-measures-new}
%%%%%%%%%%%%%%%%%%%%%%%%%%%%%%%%%%%%%%%%%%%%%%%%%%

\begin{remark}[Computing error \rev{measure} estimates]
  Observing the definitions of the error \rev{measure} estimates
  \eqref{eq:error-estimate}--\eqref{eq:discretization-error-estimate}
  yields that not only the energy $
  \densintene_\arc^\level(\spacepointsset_0; \spacediff_i)$, as a part
  of the solution~$\solution$, is required to compute the estimates
  but also the values
  $ \densintene_\arc^1(\spacepointsset_0; \spacediff_i) $ and $
  \densintene_\arc^\level(\spacepointsset_0; \spacediff_{i-1}) $,
  which are not given in terms of the solution~$ \solution $. One could compute
  these values by recomputing the \eqref{eq:overall-problem} with
  appropriately modified pipe levels and step sizes. However, this is
  computationally very costly. An alternative approach is to
  explicitly solve the modified (\wrt appropriately modified model
    levels and step sizes) energy equations of the
  Systems~\eqref{eq:discretized1}--\eqref{eq:discretized3} by means of
  implicit numerical integration. Fortunately, this is not required in
  this work since the energy equations of the
  Systems~\eqref{eq:discretized1}--\eqref{eq:discretized3} together
  with Equation~\eqref{eq:state-temperature-discr}
  allow for solving them algebraically for the energies $
  \densintene_\arc^k $, $ k = 0, 1, \ldots, n $ in linear time.
\end{remark}
In the following section we present the algorithm that adaptively switches the
previously introduced models and their discretizations by means of a
switching strategy.

%%% Local Variables:
%%% mode: latex
%%% TeX-master: "dhn-adaptive"
%%% End:

%% file: proof-error-measures-new.tex
In the following, we show that the relation
\eqref{eq:dotle-for-estimator} holds for $ \spacediff_{i} \to 0 $. In
particular, we need to show that $ \errorexd_\arc(\solution) \dotle
\errorestd_\arc(\solution) $ and $ \errorexm_\arc(\solution) \dotle
\errorestm_\arc(\solution) $ hold,
% For this, we require three preliminaries.
where the relation $ f_1(x) \dotle f_2(x) $ states
that a function~$ f_2 $ is a first-order upper bound of the function $
f_1 $ if and only if $ f_1(x) \le f_2(x) + \phi(x) $ for $ x \to 0 $ and any
function $ \phi \in o(\norm[\infty]{f_2})$.
The use of first-order error bounds that are obtained by omitting
higher-order terms is standard practice in adaptive refinement
methods; see, e.g., \cite{Ver13}. In many instances one can also
obtain exact upper bounds \cite{KonGMP03}, but these are typically
far too pessimistic to be of practical use.
%\MS{@Volker: Make a comment here that we (as usual) driop %higher-order
%terms everywhere as it is standard in adaptive refinement methods}%
% Second, \Wlog, and in
% order to keep the following derivation concise, we omit the evaluation
% argument of $ \spacepointsset_0 $ of the internal energy density $
% \densintene_\arc $. This means that we have for the exact solution $
% \densintene_a^\level(\spacepointsset_0) \to \densintene_a^\level $ and
% for its approximation $ \densintene_\arc^\level(\spacepointsset_0;
% \spacediff_i) \to \densintene_a^\level(\spacediff_{i}) $. Lastly, we
% assume that $ \lim_{\spacediff_{i} \to 0}
% \densintene_\arc^\level(\spacediff_{i} ) = \densintene_\arc^\level
% \eqqcolon \densintene_\arc^\level(0) $ holds.

We first proceed by showing that $ \errorexd_\arc(\solution) \dotle
\errorestd_\arc(\solution) $ holds for $ \spacediff_{i} \to 0 $. Since
we utilize the implicit mid-point rule to get
Systems~\eqref{eq:discretized1}--\eqref{eq:discretized3} and the fact
that its discretization error is of convergence order 2 (see, \eg,
\cite{quarteroni2010numerical}) we can write that
\begin{align}
  \densintene_\arc^\level(x_k) - \densintene_\arc^\level(x_k;
  \spacediff_i) &= c^\level(x_k) \spacediff_i^2 + \mathcal{O}(\spacediff_i^3)
                  \,, \label{eq:discretization-error-convergence-dx}
  \\
  \densintene_\arc^\level(x_k) - \densintene_\arc^\level(x_k;
  \spacediff_{i-1}) &= 4 c^\level(x_k) \spacediff_i^2 + \mathcal{O}(\spacediff_i^3) \,,
  \label{eq:discretization-error-convergence-2dx}
\end{align}
where we use $ \spacediff_{i-1} = 2 \spacediff_i $. Here, the function
$ c^\level(x) $ that arises from the Taylor series expansion of the local
discretization error is \rev{independent} of $ \spacediff_i $; see, \eg,
\cite{Stoer:2002}. Computing the
difference between \eqref{eq:discretization-error-convergence-dx} and
\eqref{eq:discretization-error-convergence-2dx} yields
\begin{equation}
  \densintene_\arc^\level(x_k; \spacediff_i) -
  \densintene_\arc^\level(x_k; \spacediff_{i-1}) = 3
  c^\level(x_k) \spacediff_i^2 + \mathcal{O}(\spacediff_i^3),
  \label{eq:discretization-error-estimate-convergence}
\end{equation}
and, thus,
\begin{equation}
  \label{eq:discretization-error-extra-term}
  c^\level(x_k) \spacediff_i^2 = \frac{\densintene_\arc^\level(x_k; \spacediff_i) -
  \densintene_\arc^\level(x_k; \spacediff_{i-1})}{3} + \mathcal{O}(\spacediff_i^3).
\end{equation}
By replacing $c^\level(x_k) \spacediff_i^2$ in
\eqref{eq:discretization-error-convergence-dx} with the result of
\eqref{eq:discretization-error-extra-term}, applying the $
\infty $-norm over $\spacepointsset_0$ on both sides, and using the
triangle inequality, we find
\begin{equation*}
  \errorexd_\arc(\solution) = \norm[\infty]{\frac{\densintene_\arc^\level(x_k; \spacediff_i) -
  \densintene_\arc^\level(x_k; \spacediff_{i-1})}{3} +
  \mathcal{O}(\spacediff_i^3)}
  \leq \frac{1}{3}\errorestd_\arc(\solution) +
    \norm[\infty]{\mathcal{O}(\spacediff_i^3)}.
\end{equation*}
Since $\errorestd_\arc(\solution) \in \mathcal{O}(\spacediff_i^2)$ holds as
shown in \eqref{eq:discretization-error-estimate-convergence},
we get that $\errorexd_\arc(\solution) \dotle
\errorestd_\arc(\solution)$ holds for $ \spacediff_{i} \to 0$.

Finally, we show that $ \errorexm_\arc(\solution) \dotle
\errorestm_\arc(\solution) $.
The ideas are rather similar.
By applying the $\infty$-norm over $\spacepointsset_0$ and the
triangle inequality to the difference between
\eqref{eq:discretization-error-convergence-dx} with \eqref{eq:model1}
and current model level $\level \in
\{\eqref{eq:model1},\eqref{eq:model2},\eqref{eq:model3}\}$ we get
\begin{align*}
  \errorexm_\arc(\solution)
  &= \norm[\infty]{\densintene_\arc^\level(x_k;
    \spacediff_i) - \densintene_\arc^1(x_k;
    \spacediff_i) + (c^\level(x_k) -
    c^1(x_k)) \spacediff_i^2 +
    \mathcal{O}(\spacediff_i^3) }
  \\
  &\leq
    \errorestm_\arc(\solution) + \norm[\infty]{\mathcal{O}(\spacediff_i^2)}.
\end{align*}
Since $\errorestm_\arc(\solution) \in \mathcal{O}(1)$, we get that
$\errorexm_\arc(\solution) \dotle \errorestm_\arc(\solution)$ holds
for $ \spacediff_{i} \to 0$.

%% file: adaptive-algorithm.tex
\section{Adaptive Algorithm}
\label{sec:adaptive}

In this section, we present and analyze the adaptive optimization
algorithm. This algorithm is based on the work
in~\cite{Mehrmann_et_al:2018} and adapted for the district heating
network problem studied in this paper. The algorithm
iteratively solves the \eqref{eq:overall-problem} while adaptively
switching the pipe model levels and discretization step sizes to
achieve a locally optimal solution that is feasible w.r.t.\ to some
prescribed tolerance.
The adaptive switching is implemented via
marking and switching strategies, which are based on the error measures
presented in the previous section.
\if0 % just a repetition
For the purpose of computing the
error \rev{measure} estimates, the energy equations of the
Systems~\eqref{eq:discretized1}--\eqref{eq:discretized3} may have to
be solved. In contrast to \cite{Mehrmann_et_al:2018}, this is done
without requiring implicit numerical integration schemes.
\fi

Given an a-priori error \rev{measure} tolerance $\errortol > 0 $, our method aims
at computing a finite sequence of solutions of the nonlinear problem
\eqref{eq:overall-problem} in order to achieve a solution~$\solution
$ with an estimated average error \rev{measure} less or equal to $\errortol $. This
motivates the following definition.

\begin{definition}
  \label{def:feasibility}
  Let $\errortol > 0 $ be a given tolerance. The solution $\solution
  $ of the \eqref{eq:overall-problem} is called $\errortol $-feasible
  if
  \begin{equation*}
    \bar \errorest (\solution) \define \frac{1}{|\pipes|}
    \sum\limits_{\arc \in \pipes} \errorest_\arc(\solution) \le
    \errortol,
  \end{equation*}
  where $\bar \errorest (\solution) $ is called the \emph{total
    average error \rev{measure} estimate}.
\end{definition}
The remainder of this section is structured as follows. We first provide the
switching and marking strategies used by our algorithm in
Section~\ref{sec:adaptive:strategies}.
Then, we present the adaptive algorithm and prove its convergence in
Section~\ref{sec:adaptive:algorithm}.

%%%%%%%%%%%%%%%%%%%%%%%%%%%%%%%%%%%%%%%%%%%%%%%%%%%%%%%%%%%%%%%%%%%%%%%%%%%
\subsection{Switching and Marking Strategies}
\label{sec:adaptive:strategies}

In a nutshell, the overall algorithm follows the standard
principles of adaptive refinement methods: a problem is solved, an error
\rev{measure} is computed, elements (here pipes) are marked to be
refined, the refinement is carried out, and the new problem is
solved; see, e.g., \cite{NocSV09,Ver13}.
In this section, we describe both the rules that are used to carry
out the refinements and the strategies that are used to mark the
pipes to be refined.

We now define switching strategies to compute
new pipe levels $\level_\arc^\new $ and new step sizes $
\spacediff_\arc^\new $. Let $\errortol > 0$ be a tolerance and $
\errortune \ge 1 $ be a tuning parameter. First, we introduce the
model level switching rules. Consider the pipe sets
\begin{equation}
  \label{eq:pipes-set-decrease}
  \pipes^{>\errortol} \coloneqq \Defset{ \arc \in \pipes }{
  \errorestm_\arc(\solution; \level_\arc) - \errorestm_\arc(\solution;
  \level_\arc^\new) > \errortol }
\end{equation}
and
\begin{equation}
  \label{eq:pipes-set-increase}
  \pipes^{<\errortune\errortol} \coloneqq \Defset{ \arc \in \pipes }{
  \errorestm_\arc(\solution; \level_\arc^\new) -
  \errorestm_\arc(\solution; \level_\arc) < \errortune \errortol }.
\end{equation}
The set~$\pipes^{>\errortol}$ ($\pipes^{<\errortune\errortol}$) contains all
the pipes for which the new model level $ \level_\arc^\new $ decreases
(increases) the model error \rev{measure} estimate
compared to the current model level $ \level_\arc $ \wrt the
error \rev{measure} tolerance $ \errortol $. In order to switch-up the model level
($ \level_\arc^\new < \level_\arc $), we apply the rule
\begin{equation}
  \label{eq:up-switching}
  %\text{(up)} \quad
  \level_\arc^\new  =
  \begin{cases}
    \level_\arc - 1,
    & \level_\arc > 1 \ \wedge \ \errorestm_\arc(\solution;
    \level_\arc) - \errorestm_\arc(\solution; \level_\arc - 1) >
    \errortol,
    \\
    1,
    & \text{otherwise}.
  \end{cases}
\end{equation}
Similarly, for down-switching of the model level ($ \level_\arc^\new >
\level_\arc $), we apply the rule
\begin{equation}
  \label{eq:down-switching}
  %\text{(down)} \quad
  \level_\arc^\new = \min \Set{\level_\arc + 1, \level_{\text{max}}}
\end{equation}
with $\level_{\text{max}} = 3$ in our setting.
According to the rules defined in Section~\ref{sec:error:notation}, we
apply the following grid refinement and coarsening rule:
\begin{equation}
  \label{eq:grid-level-switching}
  \spacediff_\arc^\new =
  \begin{cases}
    1/2 \, \spacediff_\arc,
    & \text{for a grid refinement},
    \\ 2 \, \spacediff_\arc,
    & \text{for a grid coarsening}.
  \end{cases}
\end{equation}
Based on the switching strategies defined in
\eqref{eq:pipes-set-decrease}--\eqref{eq:grid-level-switching} we can
now present our marking strategies that decide for which pipes
we switch up or down the model level and for which pipes we refine
or coarsen the step size. Let the sets \mbox{$ \refine, \upmodel \subseteq
\pipes $} represent all pipes marked for grid refinement and model level
up-switching, respectively. Furthermore, let the sets $ \coarsen,
\downmodel \subseteq \pipes
$ represent all pipes marked for grid coarsening and model level
down-switching, respectively.
To avoid unnecessary switching we use threshold parameters
$\refinementParam$, $\upswitchingParam$, $\coarseningParam$,
$\downswitchingParam \in (0,1)$.
We determine~$ \refine $ and~$ \upmodel $ by finding the minimum subset
of pipes $ \arc \in \pipes $ such that
\begin{equation}
  \label{marking:theta-refine}
  \refinementParam \sum_{\arc \in \pipes}
  \errorestd_\arc(\solution) \le \sum_{\arc \in \refine} \errorestd_\arc(\solution)
\end{equation}
and
\begin{equation}
  \label{marking:phi-up}
  \upswitchingParam \sum_{\arc \in \pipes^{>
      \errortol}} \big(\errorestm_\arc(\solution; \level_\arc) -
  \errorestm_\arc(\solution; \level_\arc^\new) \big) \le \sum_{\arc \in
    \upmodel} \big( \errorestm_\arc(\solution; \level_\arc) -
  \errorestm_\arc(\solution; \level_\arc^\new) \big)
\end{equation}
are satisfied, where in \eqref{marking:phi-up}, the rule
in~\eqref{eq:up-switching} is applied. Similarly, in order to
determine~$ \coarsen $ and~$ \downmodel $, we have to find the maximum
subset of all pipes $ \arc \in \pipes $ such that
\begin{equation}
  \label{marking:theta-coarsen}
  \coarseningParam \sum_{\arc \in \pipes}
  \errorestd_\arc(\solution) \ge \sum_{\arc \in \coarsen}
  \errorestd_\arc(\solution)
\end{equation}
and
\begin{equation}
  \label{marking:phi-down}
  \downswitchingParam \sum_{\arc \in \pipes^{<\errortune\errortol}}
  \big(\errorestm_\arc(\solution; \level_\arc^\new) -
  \errorestm_\arc(\solution; \level_\arc) \big) \ge \sum_{\arc \in \downmodel}
  \big(\errorestm_\arc(\solution; \level_\arc^\new) -
  \errorestm_\arc(\solution; \level_\arc) \big)
\end{equation}
hold, where in \eqref{marking:phi-down}, the rule
in~\eqref{eq:down-switching} is applied.

\begin{remark}
  \label{rem:exactvsestimate}
  Note that Definition~\ref{def:feasibility} is based on the total
  error \rev{measure} estimate as introduced in the previous section. Since the total
  exact error \rev{measure} is upper bounded by the total error \rev{measure} estimate
  via~\eqref{eq:dotle-for-estimator} one also has
  that the solution $ \solution $ of the \eqref{eq:overall-problem} is $
  \errortol $-feasible \wrt\ the total average exact error \rev{measure} $ \bar
  \errorex(\solution) $, \ie, $ \bar \errorex(\solution) \le
  \errortol $ holds with
  where
  \begin{equation*}
    \bar \errorex (\solution) \define \frac{1}{|\pipes|}
    \sum\limits_{\arc \in \pipes} \errorex_\arc(\solution).
    \qedhere
  \end{equation*}
  Thus, whenever error \rev{measure} estimates are used for the switching and marking
  strategies, the exact error \rev{measures} can be used as well.
\end{remark}

As used before in Section \ref{sec:error:measures}, the first-order
approximation of the discretization error \rev{measure} estimator in $x \in
[0,\length_\arc]$ of a discretization scheme of order $\beta$ reads
$\discrEstimate_a(x)  \doteq c(x) \Delta x_\arc^\beta $, where
$c(x)$ is independent of $\Delta x_\arc$. This allows us to write
\begin{equation*}
  \errorest_\arc^{\text{d,new}}(x) = \left(\frac{\Delta
      x_\arc^{\text{new}}}{\Delta x_\arc}\right)^\beta
  \errorest_\arc^{\text{d}}(x)
\end{equation*}
for the new discretization error \rev{measure} estimator after a grid refinement or
coarsening.
Since the implicit mid-point rule is used in our case, $\beta = 2$
holds, leading to
\begin{equation}
  \label{eq:errorest-coarsening-refinement}
  \errorest_\arc^{\text{d,new}}(x) =
  \begin{cases}
    \errorest_\arc^{\text{d}}(x)/4, & \text{for a grid refinement},\\
    4 \, \errorest_\arc^{\text{d}}(x), & \text{for a grid coarsening}.
  \end{cases}
\end{equation}
This also naturally holds for the exact discretization error estimator
$\errorexd_\arc(x)$.

\rev{Alternative to the error measure estimates that we present, one could use
  Richardson extrapolation based on errors on different grids to
  generate error measure estimates.
  Since our estimates are straightforward, we have not proceeded in this
  way.}

%%%%%%%%%%%%%%%%%%%%%%%%%%%%%%%%%%%%%%%%%%%%%%%%%%%%%%%%%%%%%%%%%%%%%%%%%%%
\subsection{Adaptive Algorithm}
\label{sec:adaptive:algorithm}

In this section we present the adaptive optimization algorithm.
The algorithm is formally given in
Algorithm~\ref{alg:adaptive-algorithm} and described in the following.

The input of the algorithm comprise of a complete description
of the network, including initial and boundary conditions, the error \rev{measure}
tolerance $ \errortol > 0 $ as well as initial values for the
parameters $ \refinementParam^0, \upswitchingParam^0,
\coarseningParam^0, \downswitchingParam^0 \in (0,1) $, $ \errortune^0
\le 1 $, $ \safeguard^0 \in \mathbb N_+ $. The
output of the algorithm is an $ \errortol $-feasible solution $ \solution $ of
the nonlinear problem \eqref{eq:overall-problem} according to
Definition~\ref{def:feasibility}.

% Algorithm description.
The algorithm starts by initializing model levels and grid sizes for
each pipe. It then solves the \eqref{eq:overall-problem} for the first
time and checks for $\errortol$-feasibility. Since it is likely that after the first
iteration the feasibility check fails, the algorithm enters two nested
loops: the outer loop for down-switching and coarsening and
the inner loop for up-switching and refinement. In this description we will
also refer to the outer loop as the $ k $-loop and to the inner loop
as the $ j $-loop.

Next, the inner loop is entered and the up-switching and refinement
sets~$\upmodel$ and~$\refine$ are determined. This step is followed
by up-switching and refining of each pipe accordingly. Each $ j $-loop
finishes by re-solving the \eqref{eq:overall-problem} with the new
configuration \wrt pipe model levels and grid sizes and it
checks for feasibility. The inner loop continues until either a
feasible solution $ \solution $ is found or a maximum number of inner
loop iterations $ \safeguard^k $ is reached.

What follows in the outer loop is the computation of the coarsening
and down-switching sets $ \coarsen $ and $ \downmodel $,
respectively. This step is succeeded by updating the pipe model levels
and step sizes. Similar to the inner loop, the outer loop finishes by
re-solving the \eqref{eq:overall-problem} and checking for
feasibility.

% Adaptive algorithm.
\begin{algorithm}
  \newcommand{\tikzmark}[1]{\tikz[overlay, remember picture] \node (#1) {};}
  \label{alg:adaptive-algorithm}
  \caption{Adaptive Model and Discretization Level Control}

  % Set keywords for algorithmic-package.
  \SetKwProg{AdpOptProblem}{AdpOptProblem}{}{}
  \SetKwProg{AdpLvlCtrl}{AdpLvlCtrl}{}{}

  % Input & output definition.
  \medskip \KwIn{Network $ (V,A) $, initial and
  boundary conditions, error \rev{measure} tolerance $ \errortol > 0 $, initial
  parameters $ \refinementParam^0, \upswitchingParam^0,
  \coarseningParam^0, \downswitchingParam^0 \in (0,1) $, $ \errortune^0 \le 1
  $, $ \safeguard^0 \in \mathbb N_+ $}
  \KwOut{$ \errortol $-feasible solution $ \solution $ of
    \eqref{eq:overall-problem}} \medskip\hrule\medskip

  \ForEach{$ \arc \in \pipes $}{
    Initialize model level $ \level_\arc^0 $ and step size $
    \spacediff_\arc^0 $ }
  $ y^0 \leftarrow $ Solve \eqref{eq:overall-problem} \\
  \If{$ \solution^0 $ is $ \errortol $-feasible}{
    \KwRet{$ \solution \leftarrow \solution^0 $}
  }
  \For{$ k = 1, 2, \ldots $}{
    Update parameters $ \refinementParam^k, \upswitchingParam^k,
    \coarseningParam^k, \downswitchingParam^k, \safeguard^k,
    \errortune^k $ \\
    \For{$ j = 1,\ldots, \safeguard^k $ }{
      Compute sets $ \upmodel^{k,j}, \refine^{k,j} \subseteq \pipes $
      according to \eqref{marking:theta-refine},
      \eqref{marking:phi-up} \tikzmark{UR} \\
      \ForEach{$ \arc \in \upmodel^{k,j} $}{
        Switch-up the model level $ \level_\arc^{k,j} $ according to
        \eqref{eq:up-switching}
        \label{line:switch-up}
      }
      \ForEach{$ \arc \in \refine^{k,j} $}{
        Refine step size $ \spacediff_\arc^{k,j} $ according to
        \eqref{eq:grid-level-switching}
        \label{line:refine}
      }
      $ y^{k,j} \leftarrow $ Solve \eqref{eq:overall-problem} \\
      \If{$ \solution^{k,j} $ is $ \errortol $-feasible}{
        \KwRet{$ \solution \leftarrow \solution^{k,j} $}
      }
    }
    Compute sets $ \downmodel^{k}, \coarsen^{k} \subseteq \pipes $
    according to \eqref{marking:theta-coarsen},
    \eqref{marking:phi-down} \tikzmark{DC} \\
    \ForEach{$ \arc \in \downmodel^{k} $}{
      Switch-down the model level $ \level_\arc^k $ according to
      \eqref{eq:down-switching}
      \label{line:switch-down}
    }
    \ForEach{$ \arc \in \coarsen^{k} $}{
      Coarsen step size $ \spacediff_\arc^k $ according to
      \eqref{eq:grid-level-switching}
      \label{line:coarsen}
    }
    $ \solution^{k} \leftarrow $ Solve \eqref{eq:overall-problem}
    \begin{tikzpicture}[overlay, remember picture]
        \draw[decorate, decoration={brace, amplitude=10pt}]
        (UR) ++(.5cm, 1em) to ++(0, -2.35cm) node[xshift=25pt,
        yshift=1.175cm, rotate=90, text width=2.5cm, align=center]
        {\footnotesize up-switching \& refinement};
        \draw[decorate, decoration={brace, amplitude=10pt}]
        (DC) ++(1.55cm, 1em) to ++(0, -2.35cm) node[xshift=25pt,
        yshift=1.175cm, rotate=90, text width=2.5cm, align=center]
        {\footnotesize down-switching \& coarsening};
    \end{tikzpicture}
    \\
    \If{$ \solution^{k} $ is $ \errortol $-feasible}{
      \KwRet{$ \solution \leftarrow \solution^{k} $}
    }
  }
\end{algorithm}

We first show that the algorithm is finite if we only apply
changes to the discretization step sizes while fixing the model levels
for all pipes.

\begin{lemma}
  \label{lem:discr-decrease}
  Suppose that the model level $\level_\arc \in \{1, 2, 3\}$ is fixed for
  every pipe $\arc \in \Apipe$. Let the resulting set of model levels
  be denoted by $\mathcal{M}$. Suppose further that $\eta_\arc (y) =
  \eta^d_\arc(y) $ holds in \eqref{eq:error-estimate} and that every
  \eqref{eq:overall-problem} is solved to local optimality. Consider Algorithm
  \ref{alg:adaptive-algorithm} without applying the model switching
  steps in Lines~\ref{line:switch-up}
  and~\ref{line:switch-down}. Then, the algorithm terminates after a
  finite number of refinements in Line \ref{line:refine} and
  coarsenings in Line \ref{line:coarsen}
  with an $\varepsilon$-feasible solution w.r.t.\ model level set~$\mathcal{M}$ if
  there exists a constant $C$ > 0 such that
  \begin{equation}
    \label{eq:discretization-condition}
    \frac{1}{4} \refinementParam^k \safeguard^k \geq \coarseningParam^k + C
  \end{equation}
  holds and if the step sizes of the initial discretizations are
  chosen sufficiently small.
\end{lemma}
\begin{proof}
  We focus on the total discretization error \rev{measure} defined as
  \begin{equation*}
    \discrEstimate(y^j) \define \sum_{\arc \in \Apipe} \discrEstimate_\arc (y^j)
  \end{equation*}
  and show that this quantity is positively bounded away from zero for
  one outer-loop iteration $k$ containing $\mu$ inner refinement steps and
  one coarsening step. For the sake of simplicity we drop the $k$
  index.

  Hence, we first look at the influence of one inner refinement for-loop iteration
  $j \in \{1,\dots,\mu\}$ on $\discrEstimate(y^j)$. Thus,
  \begin{equation}
    \label{eq:refinement-one-iteration}
    \begin{split}
      & \sum_{\arc \in \Apipe} \discrEstimate_\arc (y^{j-1}) -
      \sum_{\arc \in \Apipe} \discrEstimate_\arc
      (y^{j})  \\
      = \ & \sum_{\arc \in \Apipe\setminus \refine^j} \discrEstimate_\arc
      (y^{j-1}) + \sum_{\arc \in \refine^j}
      \discrEstimate_\arc (y^{j-1})  -
      \sum_{\arc \in \Apipe \setminus \refine^j} \discrEstimate_\arc
      (y^{j}) -
      \sum_{\arc \in\refine^j} \discrEstimate_\arc (y^{j}) \\
      = \ & \sum_{\arc \in\refine^j} \discrEstimate_\arc (y^{j-1}) -
      \sum_{\arc \in\refine^j} \frac{1}{4} \discrEstimate_\arc
      (y^{j-1}) \\
      = \ & \frac{3}{4} \sum_{\arc \in\refine^j}\discrEstimate_\arc
      (y^{j-1}),
    \end{split}
  \end{equation}
  where we use that $\discrEstimate_\arc (y^{j})$
  equals $1/4$ of $\discrEstimate_\arc (y^{j-1})$ if $\Delta x_\arc$ is
  chosen small enough.

  Summing up Equation~\eqref{eq:refinement-one-iteration} over all $j \in
  \{1,\dots,\mu\} $ gives the total error \rev{measure} decrease in the inner for-loop:
  \begin{align*}
    & \sum_{j = 1}^\mu \left ( \sum_{\arc \in \Apipe} \discrEstimate_\arc (y^{j-1}) -
      \sum_{\arc \in \Apipe} \discrEstimate_\arc (y^{j}) \right) \\
    = \ &        \sum_{\arc \in
      \Apipe}
    \discrEstimate_\arc
    (y^{0})
    -
    \sum_{\arc
      \in
      \Apipe}
    \discrEstimate_\arc
    (y^{\mu})\\
    = \ &
    \frac{3}{4}
    \sum_{j
      =
      1}^\mu
    \sum_{\arc
      \in\refine^j}\discrEstimate_\arc  (y^{j-1}).
  \end{align*}
  We now focus on the final coarsening step of the outer for-loop. For
  the sake of simplicity we say that $y^{\mu+1}$ corresponds to the
  solution of the \eqref{eq:overall-problem} after the coarsening
  step.
  Thus,
  \begin{align*}
    & \sum_{\arc \in \Apipe} \discrEstimate_\arc (y^{\mu+1}) -
      \sum_{\arc \in \Apipe} \discrEstimate_\arc (y^{\mu}) \\
    = \ &  \sum_{\arc \in \Apipe \setminus \coarsen} \discrEstimate_\arc
        (y^{\mu+1}) +
        \sum_{\arc \in \coarsen} \discrEstimate_\arc (y^{\mu+1}) -
        \sum_{\arc \in \Apipe \setminus \coarsen} \discrEstimate_\arc
        (y^{\mu}) - \sum_{\arc \in \coarsen} \discrEstimate_\arc
        (y^{\mu}) \\
    = \ & 4 \sum_{\arc \in \coarsen} \discrEstimate_\arc (y^{\mu}) -
        \sum_{\arc \in \coarsen} \discrEstimate_\arc (y^{\mu})\\
     = \ & 3\sum_{\arc \in \coarsen} \discrEstimate_\arc (y^{\mu})
  \end{align*}
  holds, where we again use that $\discrEstimate_\arc( y^{\mu+1})$ equals $4
  \discrEstimate_\arc (y^{\mu})$ if $\Delta x_\arc$ is
  chosen small enough.

  We now prove that the total error \rev{measure} decrease in each iteration of the outer
  for loop of Algorithm \ref{alg:adaptive-algorithm} is positive and
  uniformly bounded away from zero. Hence, we consider
  \begin{equation*}
    \sum_{\arc \in \Apipe} \discrEstimate_\arc (y^{0})
    - \sum_{\arc \in \Apipe} \discrEstimate_\arc (y^{\mu + 1}) =
    \frac{3}{4} \sum_{j = 1}^\mu \sum_{\arc
      \in\refine^j}\discrEstimate_\arc (y^{j-1}) - 3\sum_{\arc \in \coarsen}
    \discrEstimate_\arc (y^{\mu}).
  \end{equation*}
  Then, using
  \begin{equation*}
    \discrEstimate_\arc (y^{j}) \geq \discrEstimate_\arc (y^{\mu})
    \quad \text{for all} \quad j = 1, \dots , \mu,
  \end{equation*}
  \eqref{marking:theta-refine}, \eqref{eq:discretization-condition},
  and \eqref{marking:theta-coarsen}, we obtain
  \begin{align*}
    & \frac{3}{4} \sum_{j = 1}^\mu \sum_{\arc
      \in\refine^j}\discrEstimate_\arc (y^{j-1}) \geq \frac{3}{4}
      \refinementParam \sum_{j = 1}^\mu \sum_{\arc
      \in\Apipe}\discrEstimate_\arc (y^{j-1}) \geq \frac{3}{4}
      \refinementParam \sum_{j = 1}^\mu \sum_{\arc
      \in\Apipe}\discrEstimate_\arc (y^{\mu}) \\
    = \ & \frac{3}{4} \refinementParam \mu \sum_{\arc
      \in\Apipe}\discrEstimate_\arc (y^{\mu}) \geq 3( \coarseningParam +
      C) \sum_{\arc \in\Apipe}\discrEstimate_\arc (y^{\mu}) \geq 3
      \coarseningParam \sum_{\arc \in\Apipe}\discrEstimate_\arc
      (y^{\mu}) + C |\Apipe|\errortol \\
    \geq \ & 3 \sum_{\arc \in\coarsen}\discrEstimate_\arc (y^{\mu}) + C
      |\Apipe|\errortol,
  \end{align*}
  which completes the proof.
\end{proof}

Next, we show that the algorithm is finite if we only apply
model level changes while the discretization step sizes are kept
fixed.

\begin{lemma}
  \label{lem:model-decrease}
  Suppose that the discretization stepsize $\Delta x_\arc$ is fixed for
  every pipe $\arc \in \Apipe$. Suppose further that $\errorest_\arc (y) =
  \modelEstimate_\arc(y) $ holds in \eqref{eq:error-estimate} and that every
  \eqref{eq:overall-problem} is solved to local optimality. Consider Algorithm
  \ref{alg:adaptive-algorithm} without applying the discretization
  refinements in Line~\ref{line:refine} and the coarsening step in
  Line~\ref{line:coarsen}. Then, the algorithm terminates after a
  finite number of model switches in Lines \ref{line:switch-up} and \ref{line:switch-down}
  with an $\varepsilon$-feasible solution with respect to the step sizes
  $\Delta x_\arc$, $\arc \in \Apipe$, if
  there exists a constant $C$ > 0 such that
  \begin{equation}
    \label{eq:model-condition}
    \upswitchingParam^k\safeguard^k \geq \errortune^k\downswitchingParam^k|\Apipe| + C.
  \end{equation}
\end{lemma}
The proof of this lemma is the same as in~\cite{Mehrmann_et_al:2018},
which is why we omit it here.

\begin{lemma}
  \label{lem:model-estimate-refinement-coarsening}
  Let $y^\mu$ and $y^{\mu+1}$ be the solution of the
  optimization problem before and after a refinement or
  coarsening step, respectively. Let $\discrEstimate_\arc
  (y)$ and  $\modelEstimate_\arc (y)$ be the discretization and model
  error \rev{measure} estimator for a given solution $y$ of
  \eqref{eq:overall-problem} as defined in
  \eqref{eq:discretization-error-estimate} and
  \eqref{eq:model-error-estimate}. Then, if
  \begin{equation*}
    %\label{eq:lemma-inequality-condition}
    \discrEstimate_\arc (y^\mu) \ll \modelEstimate_\arc (y^\mu)
  \end{equation*}
  is satisfied, it holds that
  \begin{equation}
    \label{eq:lemma-model-after-coarsening}
    \modelEstimate_\arc (y^{\mu+1}) = \modelEstimate_\arc (y^{\mu}).
  \end{equation}
\end{lemma}
\begin{proof}
  For $x \in \Gamma_{0}$ we introduce
  $\discrEstimate_\arc(x;\ell_\arc,\Delta x_{i})$ and
  $\modelEstimate_\arc(x;\ell_\arc,\Delta x_{i})$ as the local discretization
  error \rev{measure} estimator and the local model error \rev{measure} estimator evaluated at~$x$
  using the model level $\ell_\arc$ and the step size $\Delta x_{i}$ such
  that
  \begin{align*}
    \discrEstimate_\arc(x;\ell_\arc,\Delta x_{i})
    &\define  e^{\ell_\arc}_a(x;
      \Delta x_{i}) - e^{\ell_\arc}_a(x; \Delta x_{i-1}),
    \\
    \modelEstimate_\arc(x;\ell_\arc,\Delta x_{i})
    & \define e^{1}_a(x; \Delta x_{i}) - e^{\ell_\arc}_a(x; \Delta
      x_{i})
  \end{align*}
  holds.
  Since $\discrEstimate_\arc(x;\ell_\arc,\Delta x_{i})$ uses the same
  step sizes $\Delta x_i$ and $\Delta x_{i-1}$ for all $\ell_\arc$, we
  have
  \begin{equation}
    \label{eq:iff-discretization-ll}
    \abs{\discrEstimate_\arc(x;\ell_\arc,\Delta x_{i})} \ll
    \abs{\modelEstimate_\arc(x;\ell_\arc,\Delta x_{i})}
    \iff
    \abs{\discrEstimate_\arc(x;1,\Delta x_{i})} \ll
    \abs{\modelEstimate_\arc(x;\ell_\arc,\Delta x_{i})}.
  \end{equation}
  We now focus on the coarsening step and prove Equation
  \eqref{eq:lemma-model-after-coarsening}. The proof for the refinement step is
  analogous to the coarsening step and is therefore not presented. By
  definition and due to the coarsening step, we have
  \begin{align*}
    \modelEstimate_\arc (y^{\mu+1}) =
    & \max_{x \in \Gamma_0} \  \lvert \densintene_\arc^1(x; \spacediff_{i-1}) -
      \densintene_\arc^{\ell_\arc}(x;
      \spacediff_{i-1})  \rvert
    \\
    = \ & \max_{x \in \Gamma_0}  \  \lvert \densintene_\arc^1(x; \spacediff_{i-1}) -
          \densintene_\arc^{\ell_\arc}(x;
          \spacediff_{i-1}) +
          \densintene_\arc^1(x; \spacediff_{i})
    \\
    & \qquad \quad - \densintene_\arc^1(x;
      \spacediff_{i}) +
      \densintene_\arc^{\ell_\arc}(x; \spacediff_{i}) - \densintene_\arc^{\ell_\arc}(x;
      \spacediff_{i})  \rvert \\
    = \ & \max_{x \in \Gamma_0}  \ \lvert
          \modelEstimate_\arc(x;\ell_\arc,\Delta x_{i}) -
          \discrEstimate_\arc(x;1,\Delta x_{i}) +
          \discrEstimate_\arc(x;\ell_\arc,\Delta x_{i}) \rvert.
  \end{align*}
  Using~\eqref{eq:iff-discretization-ll}, we finally obtain
  \begin{equation*}
    \modelEstimate_\arc (y^{\mu+1})
    = \max_{x \in \Gamma_0} \ \lvert
    \modelEstimate_\arc(x;\ell_\arc,\Delta
    x_{i}) \rvert
    \enifed \modelEstimate_\arc (y^{\mu}).
    \qedhere
  \end{equation*}
\end{proof}

We also have a corresponding result for the estimators of the
discretization error \rev{measure}.
For this result, we make the following assumption.

\begin{assumption}
  \label{ass:sensitivity-assumption}
  Let $y^\mu$ and $y^{\mu+1}$ be the solution of the
  optimization problem before and after a model up- or
  down-switching step, respectively.
  Moreover, let us denote with $\lambda^\mu$ and $\lambda^{\mu+1}$ the
  corresponding sensitivities.
  Then, there exists a constant $C > 0$ with $\| \lambda^\mu -
  \lambda^{\mu+1} \| \leq C$.
\end{assumption}

Before we now state the next lemma, we brief\/ly discuss this
assumption.
Informally speaking, it states that the difference of the
sensitivities (i.e., of the dual variables) of the optimization problems
before and after a model up- or down-switching step is bounded by a
constant. We are convinced that this assumption holds for the
different models in our catalog.

\begin{lemma}
  \label{lem:refinement-estimate-model-up}
  Let $y^\mu$ and $y^{\mu+1}$ respectively be the solution of the
  optimization problem before and after a model up or
  down switching step. Let $\discrEstimate_\arc
  (y)$ and $\modelEstimate_\arc (y)$ be the discretization and model
  error \rev{measure} estimator for a given solution $y$ of
  \eqref{eq:overall-problem} as defined in
  \eqref{eq:discretization-error-estimate} and
  \eqref{eq:model-error-estimate}.
  Finally, suppose that Assumption~\ref{ass:sensitivity-assumption} holds.
  Then,
  \begin{equation}
    \discrEstimate_\arc (y^{\mu+1}) = \discrEstimate_\arc (y^{\mu})
  \end{equation}
  holds.
\end{lemma}
\begin{proof}
  As long as Assumption~\ref{ass:sensitivity-assumption} holds,
  the error \rev{measure} estimate for the discretization error \rev{measure} is independent of
  the used model and we immediately get the desired result.
\end{proof}

We are now ready to prove our main theorem on the finiteness of the
proposed algorithm.

\begin{theorem}[Finite termination]
  \label{thm:finite-termination}
  Suppose that $ \discrEstimate_\arc\ll \modelEstimate_\arc$ for every
  $\arc \in \pipes$ and
  that every \eqref{eq:overall-problem} is solved to local
  optimality.
  Moreover, suppose that Assumption~\ref{ass:sensitivity-assumption} holds.
  Then, Algorithm \ref{alg:adaptive-algorithm}
  terminates after a finite number of refinements, coarsenings, and
  model switches in Lines \ref{line:switch-up}, \ref{line:refine},
  \ref{line:switch-down}, and \ref{line:coarsen} with an
  $\varepsilon$-feasible solution w.r.t.\ the reference problem if
  there exist constants $C_1 , C_2 > 0$ such that
  \begin{equation*}
    \frac{1}{4} \refinementParam^k \safeguard^k \geq
    \coarseningParam^k + C_1 \quad \text{and} \quad
    \upswitchingParam^k\safeguard^k \geq
    \errortune^k\downswitchingParam^k|\Apipe| + C_2
  \end{equation*}
  hold for all $k$.
\end{theorem}
\begin{proof}
  We first focus on the average total error \rev{measure} estimator decrease between two
  subsequent inner loop iterations of Algorithm
  \ref{alg:adaptive-algorithm}. Hence,
  \begingroup
  \allowdisplaybreaks
  \begin{align*}
   \bar {\eta} (y^{j-1}) - \bar{\eta} (y^{j}) =& \sum_{\arc \in \pipes}
                                    \eta_\arc(y^{j-1}) - \sum_{\arc
                                    \in \pipes} \eta_\arc(y^{j}) \\
    =& \sum_{\arc \in \pipes} \modelEstimate_\arc(y^{j-1}) +
      \sum_{\arc \in \pipes} \discrEstimate_\arc(y^{j-1}) - \sum_{\arc
      \in \pipes} \modelEstimate_\arc(y^{j}) - \sum_{\arc \in \pipes}
      \discrEstimate_\arc(y^{j}) \\
    =& \sum_{\arc \in \pipes \setminus ( \refine_j \cup \upmodel_j)}
       \modelEstimate_\arc(y^{j-1})
       + \sum_{\arc \in  \upmodel_j} \modelEstimate_\arc(y^{j-1})
       + \sum_{\arc \in  \refine_j \setminus  \upmodel_j }
       \modelEstimate_\arc(y^{j-1})\\
    & \quad - \sum_{\arc \in \pipes \setminus ( \refine_j \cup \upmodel_j)}
       \modelEstimate_\arc(y^{j})
       - \sum_{\arc \in  \upmodel_j} \modelEstimate_\arc(y^{j})
       - \sum_{\arc \in  \refine_j \setminus  \upmodel_j }
       \modelEstimate_\arc(y^{j})\\
    & \quad + \sum_{\arc \in \pipes \setminus ( \refine_j \cup
      \upmodel_j)} \discrEstimate_\arc(y^{j-1})
      + \sum_{\arc \in \refine_j} \discrEstimate_\arc(y^{j-1})
      + \sum_{\arc \in \upmodel_j \setminus \refine_j
      } \discrEstimate_\arc(y^{j-1})\\
    & \quad - \sum_{\arc \in \pipes \setminus ( \refine_j \cup
      \upmodel_j)} \discrEstimate_\arc(y^{j})
      - \sum_{\arc \in \refine_j} \discrEstimate_\arc(y^{j})
      - \sum_{\arc \in \upmodel_j \setminus \refine_j}
      \discrEstimate_\arc(y^{j})\\
    =& \sum_{\arc \in  \upmodel_j} \modelEstimate_\arc(y^{j-1})
       - \sum_{\arc \in  \upmodel_j} \modelEstimate_\arc(y^{j})
       + \sum_{\arc \in \refine_j} \discrEstimate_\arc(y^{j-1})
       - \sum_{\arc \in \refine_j} \discrEstimate_\arc(y^{j})\\
    =& \sum_{\arc \in  \upmodel_j} \modelEstimate_\arc(y^{j-1})
       - \sum_{\arc \in  \upmodel_j} \modelEstimate_\arc(y^{j})
       + \sum_{\arc \in \refine_j}\frac{3}{4} \discrEstimate_\arc(y^{j-1})
  \end{align*}
  \endgroup
  holds, where we use Lemma~\ref{lem:model-estimate-refinement-coarsening},
  Lemma~\ref{lem:refinement-estimate-model-up}, and
  Equation~\eqref{eq:errorest-coarsening-refinement}.
  Taking the sum over all $j = 1, \dots, \mu$ inner loop iterations gives
  \begin{align*}
    & \sum_{j = 1}^{\mu}\bar {\eta} (y^{j-1}) - \bar{\eta} (y^{j}) \\
    = \ & \bar{\eta} (y^{0}) - \bar{\eta} (y^{\mu})
    \\
    = \ & \sum_{j = 1}^{\mu} \left ( \sum_{\arc \in  \upmodel_j}
          \modelEstimate_\arc(y^{j-1})
          - \sum_{\arc \in  \upmodel_j} \modelEstimate_\arc(y^{j})
          + \sum_{\arc \in \refine_j}\frac{3}{4}
          \discrEstimate_\arc(y^{j-1}) \right ).
  \end{align*}

  Next, we focus on the outer loop iterations of Algorithm
  \ref{alg:adaptive-algorithm}. We evaluate the average total error \rev{measure} increase due to
  the coarsening and down-switching.
  Hence,
  \begingroup
  \allowdisplaybreaks
  \begin{align*}
    \bar {\eta} (y^{\mu + 1}) - \bar{\eta} (y^{\mu}) =
    & \sum_{\arc \in \pipes} \modelEstimate_\arc(y^{\mu+1}) +
      \sum_{\arc \in \pipes} \discrEstimate_\arc(y^{\mu+1}) - \sum_{\arc
      \in \pipes} \modelEstimate_\arc(y^{\mu}) - \sum_{\arc \in
      \pipes} \discrEstimate_\arc(y^{\mu}) \\
    =& \sum_{\arc \in \pipes \setminus ( \coarsen \cup \downmodel)}
       \modelEstimate_\arc(y^{\mu + 1})
       + \sum_{\arc \in  \downmodel} \modelEstimate_\arc(y^{\mu + 1})
       + \sum_{\arc \in  \coarsen \setminus  \downmodel }
       \modelEstimate_\arc(y^{\mu + 1})\\
    & \quad - \sum_{\arc \in \pipes \setminus ( \coarsen \cup \downmodel)}
       \modelEstimate_\arc(y^{\mu})
       - \sum_{\arc \in  \downmodel} \modelEstimate_\arc(y^{\mu})
       - \sum_{\arc \in  \coarsen \setminus  \downmodel }
      \modelEstimate_\arc(y^{\mu})\\
    & \quad + \sum_{\arc \in \pipes \setminus ( \coarsen \cup \downmodel)}
       \discrEstimate_\arc(y^{\mu + 1})
       + \sum_{\arc \in  \coarsen} \discrEstimate_\arc(y^{\mu + 1})
       + \sum_{\arc \in  \downmodel \setminus  \coarsen }
       \discrEstimate_\arc(y^{\mu + 1})\\
    & \quad - \sum_{\arc \in \pipes \setminus ( \coarsen \cup \downmodel)}
       \discrEstimate_\arc(y^{\mu})
       - \sum_{\arc \in  \coarsen} \discrEstimate_\arc(y^{\mu})
       - \sum_{\arc \in  \downmodel \setminus  \coarsen }
      \discrEstimate_\arc(y^{\mu})\\
    = & \sum_{\arc \in  \downmodel} \modelEstimate_\arc(y^{\mu + 1}) -
        \sum_{\arc \in  \downmodel} \modelEstimate_\arc(y^{\mu}) +
        \sum_{\arc \in  \coarsen} \discrEstimate_\arc(y^{\mu + 1}) -
        \sum_{\arc \in  \coarsen} \discrEstimate_\arc(y^{\mu}) \\
    = & \sum_{\arc \in  \downmodel} \modelEstimate_\arc(y^{\mu + 1}) -
        \sum_{\arc \in  \downmodel} \modelEstimate_\arc(y^{\mu}) +
        3 \sum_{\arc \in  \coarsen} \discrEstimate_\arc(y^{\mu}),
  \end{align*}
  \endgroup
  where we use Lemma~\ref{lem:model-estimate-refinement-coarsening},
  Lemma~\ref{lem:refinement-estimate-model-up}, and Equation
  \eqref{eq:errorest-coarsening-refinement}.

  It suffices to prove that the inner loop average total error \rev{measure} decrease is always
  greater than the outer loop average total error \rev{measure} increase, \ie,
  \begin{align*}
    & \sum_{j = 1}^{\mu} \left ( \sum_{\arc \in  \upmodel_j} \modelEstimate_\arc(y^{j-1})
       - \sum_{\arc \in  \upmodel_j} \modelEstimate_\arc(y^{j})
       + \sum_{\arc \in \refine_j}\frac{3}{4}
       \discrEstimate_\arc(y^{j-1}) \right ) \\
    > \ & \sum_{\arc \in  \downmodel} \modelEstimate_\arc(y^{\mu + 1}) -
        \sum_{\arc \in  \downmodel} \modelEstimate_\arc(y^{\mu}) +
        3 \sum_{\arc \in  \coarsen} \discrEstimate_\arc(y^{\mu}).
  \end{align*}
  Using the proofs of Lemma~\ref{lem:discr-decrease}
  and~\ref{lem:model-decrease}, we obtain
  \begin{align*}
    & \sum_{j = 1}^{\mu} \left ( \sum_{\arc \in  \upmodel_j} \modelEstimate_\arc(y^{j-1})
       - \sum_{\arc \in  \upmodel_j} \modelEstimate_\arc(y^{j})
       + \sum_{\arc \in \refine_j}\frac{3}{4}
      \discrEstimate_\arc(y^{j-1}) \right ) \\
    \geq \ & \upswitchingParam
             \mu \varepsilon + \frac{3}{4} \refinementParam \mu \sum_{\arc
             \in\Apipe}\discrEstimate_\arc (y^{\mu}) \\
    \geq \ & (\tau \downswitchingParam |\pipes| + C_2)\varepsilon
             + 3 (\coarseningParam + C_1) \sum_{\arc \in\Apipe}\discrEstimate_\arc
             (y^{\mu}) \\
    \geq \ &\sum_{\arc \in  \downmodel} (\modelEstimate_\arc(y^{\mu + 1}) -
           \modelEstimate_\arc(y^{\mu})) + C_2 \varepsilon
             + 3 \sum_{\arc \in
             \coarsen}\discrEstimate_\arc (y^{\mu}) + C_1
             |\Apipe|\errortol.
  \end{align*}
  This concludes the proof.
\end{proof}

%%% Local Variables:
%%% mode: latex
%%% TeX-master: "dhn-adaptive"
%%% End:

%% file: numerical-results.tex
\section{Numerical Results}
\label{sec:numerical-results}

In this section we present numerical results and for this we first
discuss the software and hardware setup.
Then, the considered instances are presented and, afterward,
the parameterization of the adaptive algorithm is explained.
%Finally, we present the numerical results.

%%%%%%%%%%%%%%%%%%%%%%%%%%%%%%%%%%%%%%%%%%%%%%%%%%%%%%%%%%%%%%%%%%%%%%%%%%%
\subsection{Software and Hardware Setup}
\label{sec:software-hardware}

We implemented the models in \Python~3.7.4 using the \Pyomo~6.2
package \cite{hart2011pyomo,hart2017pyomo} and solve the
resulting NLPs using the NLP solver
\CONOPTfour~4.24 \cite{drud1994conopt}, which is interfaced via the
\Pyomo-\GAMS interface.
We also tested other solvers and concluded that \CONOPTfour is the
most reliable solver that performs best for our application.
We used the default \GAMS settings.
The computations were executed on a computer with an Intel(R) Core(TM)
i7-8550U processor with eight threads at \SI{1.90}{\giga\hertz} and
\SI{16}{\giga\byte}~RAM.

%%%%%%%%%%%%%%%%%%%%%%%%%%%%%%%%%%%%%%%%%%%%%%%%%%%%%%%%%%%%%%%%%%%%%%%%%%%
\subsection{Test Instances}
\label{sec:test-case}

\begin{table}
  \centering
  \caption{Characteristics of the test networks.}
  \label{tab:network-data}
  \begin{tabular}{ccccc}
    \toprule
    Network & \# pipes & \# depots & \# consumers & total pipe length
                                                    (\si{\meter})\\
    \midrule
    \AROMA & 18 & 1 & 5 & \num{7262.4 }\\
    \STREET & 162 & 1 & 32 & \num{7627.1}\\
    \bottomrule
  \end{tabular}
\end{table}

The two networks considered in this section are the so-called
\AROMA and \STREET networks; see also \cite{Krug_et_al:2019} where they
have been used as well.
\AROMA is an academic test network, whereas \STREET is a part of an
existing real-world district heating network.
Both networks contain cycles but the much larger \STREET network only
contains a single cycle so that the overall network is almost
tree-shaped.
Table~\ref{tab:network-data} shows the main characteristics of
these networks.

The cost of waste incineration, of natural gas, and of increasing the
pressure of the water in the depot are taken from \cite{NUSSBAUMER2016496} and are
set to $\wastecost = \SI{0}{\text{\euro} \per \kWh}$, $\gascost =
\SI{0.0415}{\text{\euro} \per \kWh}$, and $\pressurecost =
\SI{0.165}{\text{\euro} \per \kWh}$. Additionally, the gas and
pressure power variables $\Powergas$ and $\Powerpress$ are
left unbounded above, whereas the waste power variable $\Powerwaste$
is bounded above by \SI{10}{\kilo \watt}.
Scarce waste incineration
power $\Powerwaste$ implies an increased consumption of costly power
($\Powerpress$ and $\Powergas$) to satisfy the total
customer demand and thus yields a non-trivial
optimization problem.

%%%%%%%%%%%%%%%%%%%%%%%%%%%%%%%%%%%%%%%%%%%%%%%%%%%%%%%%%%%%%%%%%%%%%%%%%%%
\subsection{Parameterization of the Algorithm}
\label{sec:algorithm-parameters}

\begin{table}
  \centering
  \caption{Parameters used for the numerical results.}
  \begin{tabular}{cr}
    \toprule
    Parameter & Value\\
    \midrule
    $\varepsilon$ & \SI{e-6}{\giga \joule \per \cubic\metre}\\
    $\refinementParam$ & \num{0.9}\\
    $\upswitchingParam$ & \num{0.4}\\
    $\coarseningParam$ & \num{0.45}\\
    $\downswitchingParam$ & \num{0.2}\\
    $\tau$ & \num{5}\\
    $\mu$ & \num{4}\\
    \bottomrule
  \end{tabular}
  \label{tab:numerical-parameters}
\end{table}

Table~\ref{tab:numerical-parameters} shows the parameters used for
obtaining the numerical results.
These parameters are kept constant over the course of the
iterations of the algorithm to simplify the interpretation of the
results.
It should be noted that the parameters do not satisfy the second
inequality in Theorem~\ref{thm:finite-termination}.
We choose this parameterization despite this fact because the
algorithm still converges using these settings and allows for
switching down the model level of more pipes and, hence, keeps the
optimization model more tractable over the course of the iterations.
One could, e.g., by increasing $\mu$, easily satisfy both inequalities
of Theorem~\ref{thm:finite-termination}.
For the first iteration of the adaptive
algorithm we use $\spacediff_\arc = \length_\arc/2$ and $\ell_\arc =
3$ for all $\arc \in \pipes$. This forces us to take
the reference grid~$\spacepointsset_0 = \{0,\length_\arc\}$ for all
$\arc \in \pipes$.
The assumption that the initial granularity of the discretization is
sufficiently fine is not satisfied here but does (in practice) not
harm the overall convergence of the algorithm and is therefore kept
large.

%%%%%%%%%%%%%%%%%%%%%%%%%%%%%%%%%%%%%%%%%%%%%%%%%%%%%%%%%%%%%%%%%%%%%%%%%%%
\subsection{Discussion of the Results Obtained by Using Error \rev{Measure} Estimators}
\label{sec:results-discussion-estimators}

Let us first note that none of the tested optimization solvers
converges to a feasible point for both the \AROMA and the \STREET
network when using \eqref{eq:model1} and $\spacediff_\arc
= \length_\arc / 10$ for all $\arc \in \pipes$ since this spatial
discretization already leads to a highly nonlinear problem of a size
that is very hard to be tackled by state-of-the-art NLP solvers.

\begin{figure}
  \centering
  \resizebox{0.49\textwidth}{!}{\input{figures/pgf/error-estimator/Aroma_error_values.pgf}}
  \resizebox{0.49\textwidth}{!}{\input{figures/pgf/error-estimator/street_error_values.pgf}}
  \resizebox{0.49\textwidth}{!}{\input{figures/pgf/error-estimator/Aroma_solving_times.pgf}}
  \resizebox{0.49\textwidth}{!}{\input{figures/pgf/error-estimator/street_solving_times.pgf}}
  \caption{Error estimator values (top) and computation times
    (bottom) over the course of the iterations of the adaptive
    algorithm using error \rev{measure} estimators; \AROMA network (left) and
    \STREET network (right).}
  \label{fig:estimator-error-time}
\end{figure}
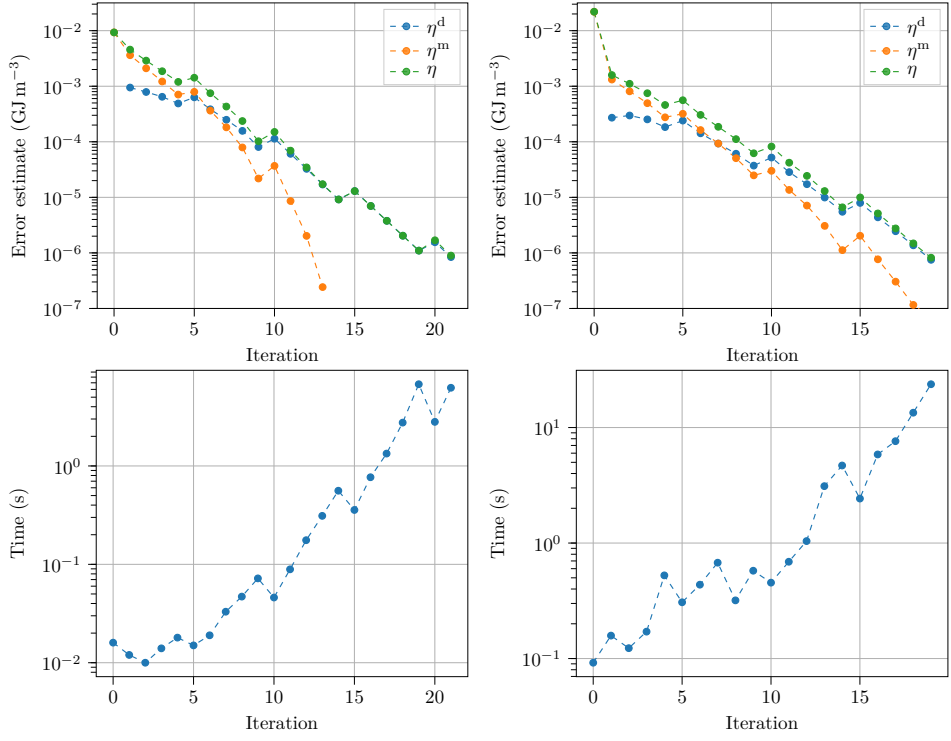

The two upper plots in Figure~\ref{fig:estimator-error-time} show
a steady decrease of the values of the error \rev{measure} estimators over the course of
the iterations of the adaptive algorithm.
Small increases in the error \rev{measure} can be observed every five iterations of
the algorithm.
These arise from the increase  of the model level and the coarsening of the
discretizations (outer loop) that is carried out after
four refinement steps in which we increase the model's accuracy (inner
loop).
The error \rev{measure} plots thus confirm that the algorithm steadily
decreases the total error \rev{measure} over the course of one outer loop
iteration.
\if0 %%%%%%%%%%%%%%%%%%%%%%%%%%%%%%%%%%%%%%%%%%%%%%%%%%
Moreover, for both networks the average convergence rate of the model
error \rev{measure} estimators increases for larger iterations whereas the
convergence rate of the discretization error \rev{measure} estimators
appears to be constant. This confirms the aspects discussed
in Remark~\ref{remark:convergence-rate-error-estimators}.}
\fi %%%%%%%%%%%%%%%%%%%%%%%%%%%%%%%%%%%%%%%%%%%%%%%%%%

The results show that the algorithm works as expected and that it
terminates after a finite number of iterations with a locally optimal
solution of a model that
has a physical accuracy for which state-of-the-art solvers are not
able to compute a solution from scratch.
This is one of the most important contributions of this paper:
We can solve realistic instances that have not been solvable before.
Additionally, the two lower plots in
Figure~\ref{fig:estimator-error-time} show the computation times for
the separate models of Type~\eqref{eq:overall-problem} that we solve
in every iteration.
Although we warmstart every new problem with the solution of the
previous one, we observe an increase of solution times due to the
higher complexity of the successive models that we solve.

\begin{figure}
  \centering
  \resizebox{0.49\textwidth}{!}{\input{figures/pgf/error-estimator/Aroma_inner_outer_loop_sets.pgf}}
  \resizebox{0.49\textwidth}{!}{\input{figures/pgf/error-estimator/street_inner_outer_loop_sets.pgf}}
  \caption{Proportion of pipes inside sets
    $\upmodel,\refine,\downmodel$, and $\coarsen$ over the course of
    iterations of the adaptive algorithm using the error \rev{measure} estimators;
    \AROMA (left) and \STREET (right).}
  \label{fig:estimator-sets}
\end{figure}

\begin{figure}
  \centering
  \resizebox{0.49\textwidth}{!}{\input{figures/pgf/error-estimator/street_model_level_sets_sizes.pgf}}
  \resizebox{0.49\textwidth}{!}{\input{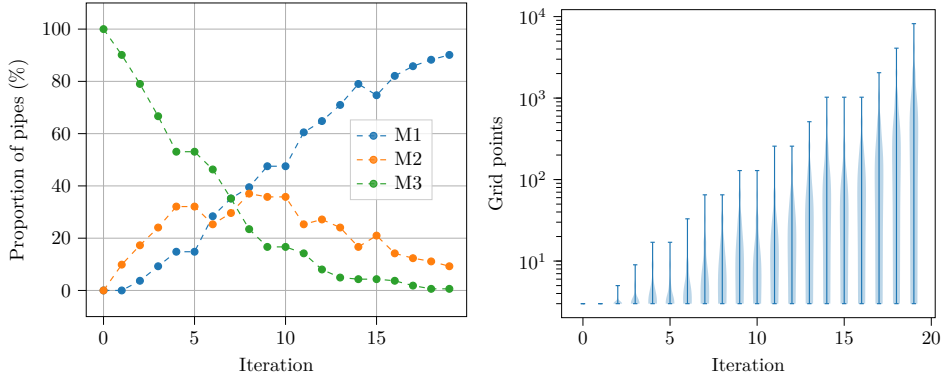}}
  \caption{Proportion of pipes inside each model level set (left) and violin
    plots of the quantity of grid points (right) over
    the iterations of our adaptive algorithm  using the error \rev{measure} estimators
    applied on the \STREET network.}
  \label{fig:estimator-level-discr}
\end{figure}

Next, Figure~\ref{fig:estimator-sets} shows the proportion of pipes
inside the sets~$\upmodel$, $\refine$, $\downmodel$, and $\coarsen$
before solving \eqref{eq:overall-problem} for every iteration of the
algorithm.
The discretization sets represent a larger proportion of pipes when
compared to the sets for switching between the model levels.
This originates from the parameter selection that favors
changes of the discretization and is explained by the fact that the model
level of a specific pipe can only be increased twice---unlike the
discretization step size that may need to be halved more often.
The right plot of Figure~\ref{fig:estimator-level-discr}
shows violin plots for the amount of grid points in the pipes over the
iterations of the algorithm applied to the \STREET network.
The plot confirms the idea behind the parameter selection.
Besides this, Figure~\ref{fig:estimator-sets} illustrates that the
down-switching set $\downmodel$ stays empty until the last outer loop
iteration for both networks.
This is a result of the set $\pipes^{<\errortune\errortol}$ being
empty for the first outer-loop iterations of the algorithm, which
forces $\downmodel$ to be empty.
The amount of pipes in each model level is shown in the left plot of
Figure~\ref{fig:estimator-level-discr}.
Roughly \SI{90}{\percent} of all pipes end up in the most accurate
model level whereas the remaining stay in the intermediate level.

Overall, we see that the behavior of the algorithm is comparable when
applied to the two different networks, which indicates that the
algorithm is robust.

%%%%%%%%%%%%%%%%%%%%%%%%%%%%%%%%%%%%%%%%%%%%%%%%%%%%%%%%%%%%%%%%%%%%%%%%%%%
\subsection{Discussion of the Results Obtained by Using Exact Error \rev{Measures}}
\label{sec:results-discussion-exact-errors}

\begin{figure}
  \centering
  \resizebox{0.49\textwidth}{!}{\input{figures/pgf/exact-error/street_error_values.pgf}}
  \resizebox{0.49\textwidth}{!}{\input{figures/pgf/exact-error/street_solving_times.pgf}}
  \resizebox{0.49\textwidth}{!}{\input{figures/pgf/exact-error/street_model_level_sets_sizes.pgf}}
  \resizebox{0.49\textwidth}{!}{\input{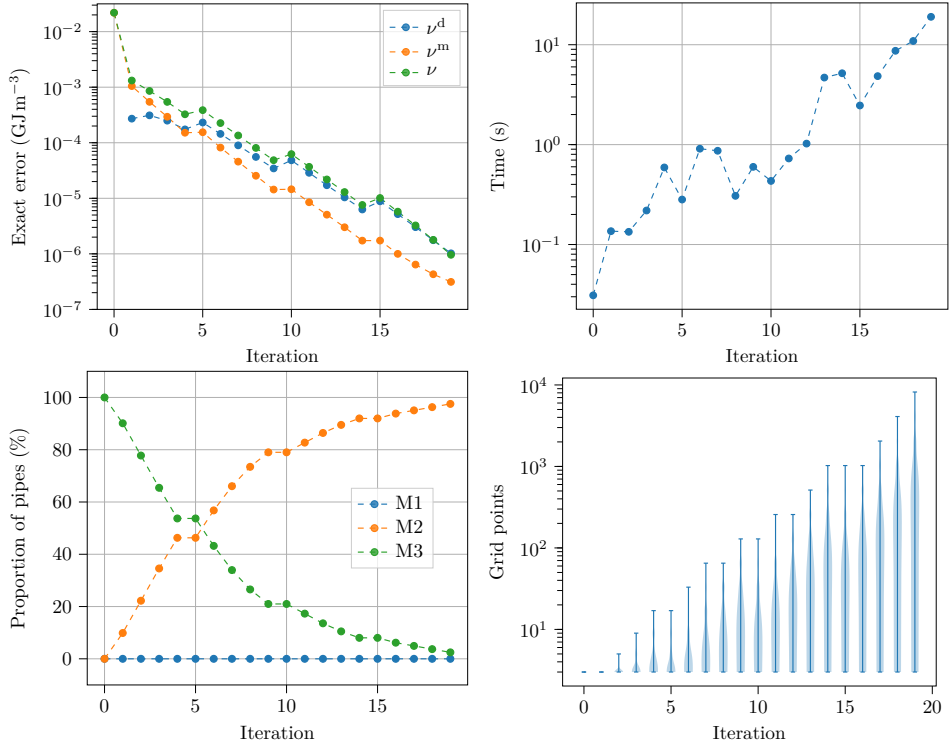}}
  \caption{Exact error values (top-left), computation time
    (top-right), proportion of pipes inside each model level set
    (bottom-left), and violin plots of the quantity of grid points
    (bottom-right) over the course of the iterations of the adaptive
    algorithm using the exact error \rev{measures} applied on the \STREET network.}
  \label{fig:exact-error-comparison}
\end{figure}

We now compare the impact of using the error \rev{measure} estimators defined in
\eqref{eq:error-estimate}--\eqref{eq:discretization-error-estimate}
when employing the exact error \rev{measures} as defined in
\eqref{eq:exact-total-error}--\eqref{eq:exact-discretization-error}.
To this end, we only consider the larger \STREET network.
Figure~\ref{fig:exact-error-comparison} shows the previously discussed
plots using exact error \rev{measures}.
Both approaches need 19~iterations to reach the desired tolerance.
However, when looking at the distribution of model levels, we see that
in the case of using error \rev{measure} estimators, a much higher proportion of
pipes are modeled using the most accurate model~\eqref{eq:model1},
which is not the case for any pipe in the exact error \rev{measure} case; see the
bottom-left plot in Figure~\ref{fig:exact-error-comparison}.
Thus, it seems that the error \rev{measure} estimators overestimate the importance
of switching to the most accurate model level.
Consequently, using the error \rev{measure} estimators instead of the exact error \rev{measures} introduces a larger amount of nonlinearities to the models that are
solved in each iteration.
This is an interesting aspect and shows that it might be beneficial to
use exact error \rev{measures} if they are available like for the ODEs that we consider
in this paper.
Nevertheless, the computation times show very similar behavior for
both approaches, which makes clear that using error \rev{measure} estimators
(especially in cases in which exact error \rev{measure} formulas are not available)
also leads to an effective method.

%%%%%%%%%%%%%%%%%%%%%%%%%%%%%%%%%%%%%%%%%%%%%%%%%%%%%%%%%%%%%%%%%%%%%%%%%%%
\subsection{Is Physical Accuracy Worth the Effort?}
\label{sec:phys-accur-worth}

\begin{figure}
  \centering
  % \resizebox{\textwidth}{!}{\input{figures/tikz/aroma_visu_first.tikz}} % tikz file
  \resizebox{\textwidth}{!}{\includegraphics{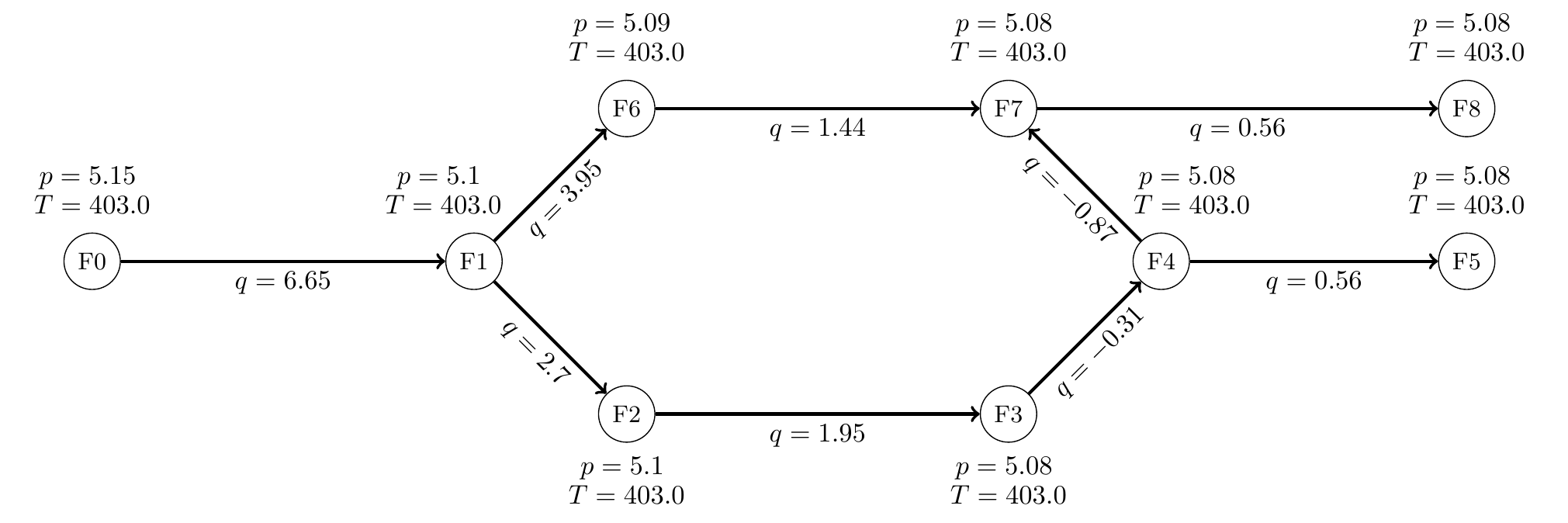}} % "tikz file converted to pdf" fallback

  \vspace*{1em}

  % \resizebox{\textwidth}{!}{\input{figures/tikz/aroma_visu_last.tikz}} % tikz file
  \resizebox{\textwidth}{!}{\includegraphics{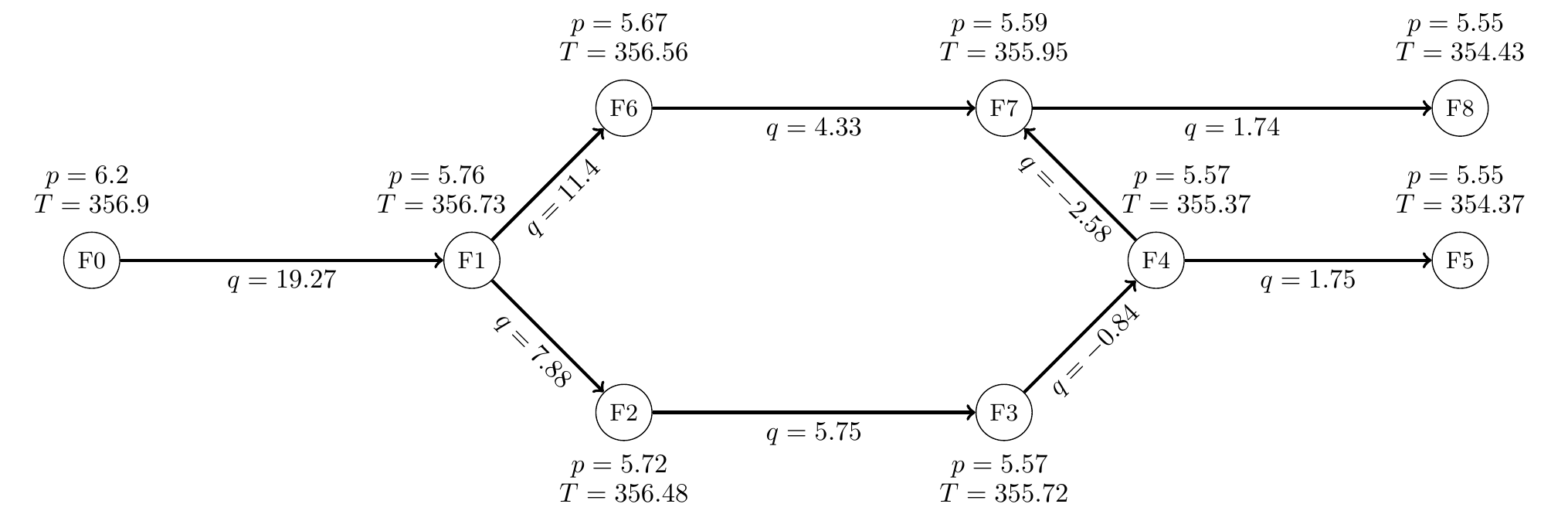}} % "tikz file converted to pdf" fallback
  \caption{Forward flow solution of the \AROMA network using the
    error \rev{measure} estimators after the first (top figure) and last (bottom
    figure) iteration of the adaptive algorithm. Temperatures in
    \si{\kelvin}, pressures in \si{\bar}, and mass flows in
    \si{\kilogram \per \second}.}
  \label{fig:aroma-first-last}
\end{figure}

Let us close this section with a brief analysis of whether the
physical accuracy guaranteed by our adaptive method is worth the
computational effort.
The answer is a clear ``Yes''.
To illustrate this, Figure~\ref{fig:aroma-first-last} shows the values
of some forward flow variables (pressures,
temperatures, and mass flows) that are part of the
\eqref{eq:overall-problem} of the \AROMA network solved in the first
iteration as well as in the last iteration of the adaptive
algorithm.
The parameter setup used in this test case is the same as presented in
Section~\ref{sec:test-case}.
Note that the solution of the first iteration (top figure)
corresponds to a rather coarse physical modeling whereas the solution
of the last iteration (bottom figure) satisfies the prescribed
tolerance and is very accurate.

The difference of the solution values are obvious.
The first solution has no temperature losses at all (see
\eqref{eq:model3}) and all temperature values are at the upper bound.
Moreover, the mass flow values are comparably small.
This changes completely in the final solution.
The temperatures have decreased around \SI{50}{\kelvin} and mass flows
have increased by up to a factor of~3.
The pressures have also changed by around \SI{10}{\percent}.
It is clearly visible that the physical solution and the control of
the network changes significantly if the physical accuracy is
increased.
Thus, there is a strong need for computing highly accurate solutions
if the resulting controls shall be practically useful.

%%% Local Variables:
%%% mode: latex
%%% TeX-master: "dhn-adaptive"
%%% End:

%% file: figures/pgf/error-estimator/Aroma_error_values.pgf
% This file was created with tikzplotlib v0.9.15.
\begin{tikzpicture}

\definecolor{color0}{rgb}{0.12156862745098,0.466666666666667,0.705882352941177}
\definecolor{color1}{rgb}{1,0.498039215686275,0.0549019607843137}
\definecolor{color2}{rgb}{0.172549019607843,0.627450980392157,0.172549019607843}

\begin{axis}[
legend cell align={left},
legend style={fill opacity=0.8, draw opacity=1, text opacity=1, draw=white!80!black},
log basis y={10},
tick align=outside,
tick pos=left,
x grid style={white!69.0196078431373!black},
xlabel={Iteration},
xmajorgrids,
xmin=-1.05, xmax=22.05,
xtick style={color=black},
y grid style={white!69.0196078431373!black},
ylabel={Error estimate (\si{\giga\joule\per\cubic\metre})},
ymajorgrids,
ymin=1e-07, ymax=0.0316227766016838,
ymode=log,
ytick style={color=black}
]
\addplot [semithick, color0, dashed, mark=*, mark size=1.75, mark options={solid}]
table {%
0 0
1 0.00094698825649238
2 0.000787983168905442
3 0.000645626325391885
4 0.000489640995590238
5 0.000631171515724541
6 0.000385040023112125
7 0.000249012801074502
8 0.000157309126480021
9 8.0521396609106e-05
10 0.000113695460610681
11 6.06940095134411e-05
12 3.25712955196086e-05
13 1.70172338676088e-05
14 9.15779349156055e-06
15 1.29701237685929e-05
16 6.99934244735705e-06
17 3.78179209086299e-06
18 2.03862992738187e-06
19 1.09902727365494e-06
20 1.55400266064372e-06
21 8.41320828187797e-07
};
\addlegendentry{$\eta^{\text{d}}$}
\addplot [semithick, color1, dashed, mark=*, mark size=1.75, mark options={solid}]
table {%
0 0.00930866594511284
1 0.00361418818966633
2 0.00210324913354764
3 0.0012171898397333
4 0.000707640965034149
5 0.000792087138720304
6 0.000363680169234811
7 0.000181932467528802
8 7.90894697785891e-05
9 2.18318505259107e-05
10 3.68652215642333e-05
11 8.59447768468824e-06
12 2.02360122703844e-06
13 2.42412119426661e-07
};
\addlegendentry{$\eta^{\text{m}}$}
\addplot [semithick, color2, dashed, mark=*, mark size=1.75, mark options={solid}]
table {%
0 0.00930866594511284
1 0.00456117644615871
2 0.00289123230245308
3 0.00186281616512519
4 0.00119728196062439
5 0.00142325865444485
6 0.000748720192346936
7 0.000430945268603303
8 0.000236398596258611
9 0.000102353247135017
10 0.000150560682174914
11 6.92884871981293e-05
12 3.4594896746647e-05
13 1.72596459870355e-05
14 9.15779349156055e-06
15 1.29701237685929e-05
16 6.99934244735705e-06
17 3.78179209086299e-06
18 2.03862992738187e-06
19 1.09902727365494e-06
20 1.6885394545214e-06
21 8.9483398158186e-07
};
\addlegendentry{$\eta$}
\end{axis}

\end{tikzpicture}

%% file: figures/pgf/error-estimator/street_error_values.pgf
% This file was created with tikzplotlib v0.9.15.
\begin{tikzpicture}

\definecolor{color0}{rgb}{0.12156862745098,0.466666666666667,0.705882352941177}
\definecolor{color1}{rgb}{1,0.498039215686275,0.0549019607843137}
\definecolor{color2}{rgb}{0.172549019607843,0.627450980392157,0.172549019607843}

\begin{axis}[
legend cell align={left},
legend style={fill opacity=0.8, draw opacity=1, text opacity=1, draw=white!80!black},
log basis y={10},
tick align=outside,
tick pos=left,
x grid style={white!69.0196078431373!black},
xlabel={Iteration},
xmajorgrids,
xmin=-0.95, xmax=19.95,
xtick style={color=black},
y grid style={white!69.0196078431373!black},
ylabel={Error estimate (\si{\giga\joule\per\cubic\metre})},
ymajorgrids,
ymin=1e-07, ymax=0.0316227766016838,
ymode=log,
ytick style={color=black}
]
\addplot [semithick, color0, dashed, mark=*, mark size=1.75, mark options={solid}]
table {%
0 0
1 0.000271426083620309
2 0.000297371814692151
3 0.000253143738234624
4 0.000183834813249692
5 0.000240749640648987
6 0.000142073779997202
7 9.3007901966655e-05
8 6.09319257691729e-05
9 3.72596780283184e-05
10 5.20418117735765e-05
11 2.84350854161746e-05
12 1.72424756474221e-05
13 9.92007914924033e-06
14 5.49555112435512e-06
15 7.93789717218758e-06
16 4.36737839547738e-06
17 2.45416865755618e-06
18 1.36830201378023e-06
19 7.50904893031459e-07
};
\addlegendentry{$\eta^{\text{d}}$}
\addplot [semithick, color1, dashed, mark=*, mark size=1.75, mark options={solid}]
table {%
0 0.0218973230866546
1 0.0013190350008404
2 0.000809593273130475
3 0.000495330006552579
4 0.000275714046762891
5 0.00031916600638478
6 0.00016265461195082
7 9.30779856357919e-05
8 5.06703074891648e-05
9 2.48850248824936e-05
10 3.00266580238808e-05
11 1.36058943468776e-05
12 7.11455233953544e-06
13 3.06768241645929e-06
14 1.12175205924997e-06
15 2.03386330390068e-06
16 7.67189671902193e-07
17 3.03846160505473e-07
18 1.15804896341062e-07
19 6.9219821051499e-08
};
\addlegendentry{$\eta^{\text{m}}$}
\addplot [semithick, color2, dashed, mark=*, mark size=1.75, mark options={solid}]
table {%
0 0.0218973230866546
1 0.0015904610844607
2 0.00110696508782263
3 0.000748473744787203
4 0.000459548860012583
5 0.000559915647033767
6 0.000304728391948022
7 0.000186085887602447
8 0.000111602233258338
9 6.2144702910812e-05
10 8.20684697974573e-05
11 4.20409797630522e-05
12 2.43570279869575e-05
13 1.29877615656996e-05
14 6.61730318360509e-06
15 9.97176047608827e-06
16 5.13456806737957e-06
17 2.75801481806165e-06
18 1.48410691012129e-06
19 8.20124714082958e-07
};
\addlegendentry{$\eta$}
\end{axis}

\end{tikzpicture}

%% file: figures/pgf/error-estimator/Aroma_solving_times.pgf
% This file was created with tikzplotlib v0.9.15.
\begin{tikzpicture}

\definecolor{color0}{rgb}{0.12156862745098,0.466666666666667,0.705882352941177}

\begin{axis}[
log basis y={10},
tick align=outside,
tick pos=left,
x grid style={white!69.0196078431373!black},
xlabel={Iteration},
xmajorgrids,
xmin=-1.05, xmax=22.05,
xtick style={color=black},
y grid style={white!69.0196078431373!black},
ylabel={Time (s)},
ymajorgrids,
ymin=0.00721870262801969, ymax=9.37963340570835,
ymode=log,
ytick style={color=black}
]
\addplot [semithick, color0, dashed, mark=*, mark size=1.75, mark options={solid}]
table {%
0 0.016000214964151
1 0.012000161223114
2 0.009999820031226
3 0.013999873772264
4 0.017999927513301
5 0.015000044368207
6 0.019000098109245
7 0.032999971881509
8 0.046999845653772
9 0.071999710053205
10 0.046000303700566
11 0.089000095613301
12 0.175999850034714
13 0.311000249348581
14 0.559999980032444
15 0.356999924406409
16 0.768000259995461
17 1.33600034750998
18 2.75700010824949
19 6.77100028842688
20 2.80499949585646
21 6.23199949041009
};
\end{axis}

\end{tikzpicture}

%% file: figures/pgf/error-estimator/street_solving_times.pgf
% This file was created with tikzplotlib v0.9.15.
\begin{tikzpicture}

\definecolor{color0}{rgb}{0.12156862745098,0.466666666666667,0.705882352941177}

\begin{axis}[
legend style={
  fill opacity=0.8,
  draw opacity=1,
  text opacity=1,
  at={(0.03,0.97)},
  anchor=north west,
  draw=white!80!black
},
log basis y={10},
tick align=outside,
tick pos=left,
x grid style={white!69.0196078431373!black},
xlabel={Iteration},
xmajorgrids,
xmin=-0.95, xmax=19.95,
xtick style={color=black},
y grid style={white!69.0196078431373!black},
ylabel={Time (s)},
ymajorgrids,
ymin=0.069698464424885, ymax=31.307081686677,
ymode=log,
ytick style={color=black}
]
\addplot [semithick, color0, dashed, mark=*, mark size=1.75, mark options={solid}, forget plot]
table {%
0 0.091999978758395
1 0.157999922521412
2 0.122999609448016
3 0.17100025434047
4 0.525000295601785
5 0.307000195607543
6 0.435999571345747
7 0.675999652594328
8 0.318999728187919
9 0.576000194996595
10 0.454000127501786
11 0.688999984413385
12 1.03999951388687
13 3.10900043696165
14 4.69299971591681
15 2.42699976079166
16 5.85000033024699
17 7.62099944986403
18 13.4440002031624
19 23.7180002499372
};
\end{axis}

\end{tikzpicture}

%% file: figures/pgf/error-estimator/Aroma_inner_outer_loop_sets.pgf
% This file was created with tikzplotlib v0.9.15.
\begin{tikzpicture}

\definecolor{color0}{rgb}{0.12156862745098,0.466666666666667,0.705882352941177}
\definecolor{color1}{rgb}{1,0.498039215686275,0.0549019607843137}
\definecolor{color2}{rgb}{0.172549019607843,0.627450980392157,0.172549019607843}
\definecolor{color3}{rgb}{0.83921568627451,0.152941176470588,0.156862745098039}

\begin{axis}[
legend cell align={left},
legend style={
  fill opacity=0.8,
  draw opacity=1,
  text opacity=1,
  at={(0.03,0.97)},
  anchor=north west,
  draw=white!80!black
},
tick align=outside,
tick pos=left,
x grid style={white!69.0196078431373!black},
xlabel={Iteration},
xmajorgrids,
xmin=-1, xmax=22,
xtick style={color=black},
y grid style={white!69.0196078431373!black},
ylabel={Proportion of pipes (\%)},
ymajorgrids,
ymin=-10, ymax=110,
ytick style={color=black}
]
\addplot [semithick, color0, dashed, mark=*, mark size=1.75, mark options={solid}]
table {%
1 16.6666666666667
2 22.2222222222222
3 22.2222222222222
4 22.2222222222222
6 27.7777777777778
7 22.2222222222222
8 16.6666666666667
9 16.6666666666667
11 11.1111111111111
12 11.1111111111111
13 5.55555555555556
14 5.55555555555556
16 0
17 0
18 0
19 0
21 5.55555555555556
};
\addlegendentry{$\mathcal{U}$}
\addplot [semithick, color1, dashed, mark=*, mark size=1.75, mark options={solid}]
table {%
1 0
2 16.6666666666667
3 27.7777777777778
4 38.8888888888889
6 55.5555555555556
7 66.6666666666667
8 77.7777777777778
9 88.8888888888889
11 88.8888888888889
12 88.8888888888889
13 88.8888888888889
14 88.8888888888889
16 88.8888888888889
17 88.8888888888889
18 88.8888888888889
19 88.8888888888889
21 88.8888888888889
};
\addlegendentry{$\mathcal{R}$}
\addplot [semithick, color2, dashed, mark=*, mark size=1.75, mark options={solid}]
table {%
5 0
10 0
15 0
20 11.1111111111111
};
\addlegendentry{$\mathcal{D}$}
\addplot [semithick, color3, dashed, mark=*, mark size=1.75, mark options={solid}]
table {%
5 22.2222222222222
10 50
15 50
20 50
};
\addlegendentry{$\mathcal{C}$}
\end{axis}

\end{tikzpicture}

%% file: figures/pgf/error-estimator/street_inner_outer_loop_sets.pgf
% This file was created with tikzplotlib v0.9.15.
\begin{tikzpicture}

\definecolor{color0}{rgb}{0.12156862745098,0.466666666666667,0.705882352941177}
\definecolor{color1}{rgb}{1,0.498039215686275,0.0549019607843137}
\definecolor{color2}{rgb}{0.172549019607843,0.627450980392157,0.172549019607843}
\definecolor{color3}{rgb}{0.83921568627451,0.152941176470588,0.156862745098039}

\begin{axis}[
legend cell align={left},
legend style={
  fill opacity=0.8,
  draw opacity=1,
  text opacity=1,
  at={(0.03,0.97)},
  anchor=north west,
  draw=white!80!black
},
tick align=outside,
tick pos=left,
x grid style={white!69.0196078431373!black},
xlabel={Iteration},
xmajorgrids,
xmin=-1, xmax=19.9,
xtick style={color=black},
y grid style={white!69.0196078431373!black},
ylabel={Proportion of pipes (\%)},
ymajorgrids,
ymin=-10, ymax=110,
ytick style={color=black}
]
\addplot [semithick, color0, dashed, mark=*, mark size=1.75, mark options={solid}]
table {%
1 9.87654320987654
2 14.8148148148148
3 17.9012345679012
4 19.1358024691358
6 20.3703703703704
7 17.9012345679012
8 16.0493827160494
9 14.8148148148148
11 15.4320987654321
12 10.4938271604938
13 9.25925925925926
14 8.64197530864197
16 8.02469135802469
17 5.55555555555556
18 3.7037037037037
19 1.85185185185185
};
\addlegendentry{$\mathcal{U}$}
\addplot [semithick, color1, dashed, mark=*, mark size=1.75, mark options={solid}]
table {%
1 0
2 8.64197530864197
3 17.9012345679012
4 29.0123456790123
6 40.7407407407407
7 46.2962962962963
8 56.1728395061728
9 65.4320987654321
11 72.2222222222222
12 74.0740740740741
13 79.0123456790123
14 81.4814814814815
16 82.7160493827161
17 83.3333333333333
18 85.1851851851852
19 85.8024691358025
};
\addlegendentry{$\mathcal{R}$}
\addplot [semithick, color2, dashed, mark=*, mark size=1.75, mark options={solid}]
table {%
5 0
10 0
15 4.32098765432099
};
\addlegendentry{$\mathcal{D}$}
\addplot [semithick, color3, dashed, mark=*, mark size=1.75, mark options={solid}]
table {%
5 17.2839506172839
10 41.358024691358
15 51.2345679012346
};
\addlegendentry{$\mathcal{C}$}
\end{axis}

\end{tikzpicture}

%% file: figures/pgf/error-estimator/street_model_level_sets_sizes.pgf
% This file was created with tikzplotlib v0.9.15.
\begin{tikzpicture}

\definecolor{color0}{rgb}{0.12156862745098,0.466666666666667,0.705882352941177}
\definecolor{color1}{rgb}{1,0.498039215686275,0.0549019607843137}
\definecolor{color2}{rgb}{0.172549019607843,0.627450980392157,0.172549019607843}

\begin{axis}[
legend cell align={left},
legend style={
  fill opacity=0.8,
  draw opacity=1,
  text opacity=1,
  at={(0.91,0.5)},
  anchor=east,
  draw=white!80!black
},
tick align=outside,
tick pos=left,
x grid style={white!69.0196078431373!black},
xlabel={Iteration},
xmajorgrids,
xmin=-0.95, xmax=19.95,
xtick style={color=black},
y grid style={white!69.0196078431373!black},
ylabel={Proportion of pipes (\%)},
ymajorgrids,
ymin=-10, ymax=110,
ytick style={color=black}
]
\addplot [semithick, color0, dashed, mark=*, mark size=1.75, mark options={solid}]
table {%
0 0
1 0
2 3.7037037037037
3 9.25925925925926
4 14.8148148148148
5 14.8148148148148
6 28.3950617283951
7 35.1851851851852
8 39.5061728395062
9 47.5308641975309
10 47.5308641975309
11 60.4938271604938
12 64.8148148148148
13 70.9876543209877
14 79.0123456790123
15 74.6913580246914
16 82.0987654320988
17 85.8024691358025
18 88.2716049382716
19 90.1234567901235
};
\addlegendentry{M1}
\addplot [semithick, color1, dashed, mark=*, mark size=1.75, mark options={solid}]
table {%
0 0
1 9.87654320987654
2 17.2839506172839
3 24.0740740740741
4 32.0987654320988
5 32.0987654320988
6 25.3086419753086
7 29.6296296296296
8 37.037037037037
9 35.8024691358025
10 35.8024691358025
11 25.3086419753086
12 27.1604938271605
13 24.0740740740741
14 16.6666666666667
15 20.9876543209877
16 14.1975308641975
17 12.3456790123457
18 11.1111111111111
19 9.25925925925926
};
\addlegendentry{M2}
\addplot [semithick, color2, dashed, mark=*, mark size=1.75, mark options={solid}]
table {%
0 100
1 90.1234567901235
2 79.0123456790123
3 66.6666666666667
4 53.0864197530864
5 53.0864197530864
6 46.2962962962963
7 35.1851851851852
8 23.4567901234568
9 16.6666666666667
10 16.6666666666667
11 14.1975308641975
12 8.02469135802469
13 4.93827160493827
14 4.32098765432099
15 4.32098765432099
16 3.7037037037037
17 1.85185185185185
18 0.617283950617284
19 0.617283950617284
};
\addlegendentry{M3}
\end{axis}

\end{tikzpicture}

%% file: figures/pgf/exact-error/street_error_values.pgf
% This file was created with tikzplotlib v0.9.15.
\begin{tikzpicture}

\definecolor{color0}{rgb}{0.12156862745098,0.466666666666667,0.705882352941177}
\definecolor{color1}{rgb}{1,0.498039215686275,0.0549019607843137}
\definecolor{color2}{rgb}{0.172549019607843,0.627450980392157,0.172549019607843}

\begin{axis}[
legend cell align={left},
legend style={fill opacity=0.8, draw opacity=1, text opacity=1, draw=white!80!black},
log basis y={10},
tick align=outside,
tick pos=left,
x grid style={white!69.0196078431373!black},
xlabel={Iteration},
xmajorgrids,
xmin=-0.95, xmax=19.95,
xtick style={color=black},
y grid style={white!69.0196078431373!black},
ylabel={Exact error (\si{\giga\joule\per\cubic\metre})},
ymajorgrids,
ymin=1e-07, ymax=0.0316227766016838,
ymode=log,
ytick style={color=black}
]
\addplot [semithick, color0, dashed, mark=*, mark size=1.75, mark options={solid}]
table {%
0 0
1 0.000271390142989099
2 0.000311395096218296
3 0.000249724719469138
4 0.000173979993821515
5 0.00023145734848569
6 0.000144298675972405
7 8.97481477216199e-05
8 5.55686858556708e-05
9 3.44991361198e-05
10 4.82863162026838e-05
11 2.87604625088946e-05
12 1.71841720004074e-05
13 1.03779580626973e-05
14 6.29348531686395e-06
15 8.82598405611994e-06
16 5.18982516986224e-06
17 3.02367983579709e-06
18 1.75964772931698e-06
19 1.02216205607355e-06
};
\addlegendentry{$\nu^{\text{d}}$}
\addplot [semithick, color1, dashed, mark=*, mark size=1.75, mark options={solid}]
table {%
0 0.0218808284247954
1 0.00104763465146154
2 0.000544495526922582
3 0.0002936787204453
4 0.000151369630468577
5 0.000154576921500533
6 8.18369702150303e-05
7 4.56723213305462e-05
8 2.54641509696893e-05
9 1.4375350549735e-05
10 1.45519107580523e-05
11 8.49280440575841e-06
12 5.07157555500593e-06
13 3.01114589921597e-06
14 1.72994029636957e-06
15 1.7325704892756e-06
16 9.99855154925842e-07
17 6.41038326165926e-07
18 4.29458722578155e-07
19 3.12852859030719e-07
};
\addlegendentry{$\nu^{\text{m}}$}
\addplot [semithick, color2, dashed, mark=*, mark size=1.75, mark options={solid}]
table {%
0 0.0218808284247954
1 0.00131902285528453
2 0.000855880898617631
3 0.000543380037321182
4 0.000325277831443068
5 0.000385965699518653
6 0.00022598811907715
7 0.000135192276699018
8 8.07318753037464e-05
9 4.85180592935886e-05
10 6.24897128445548e-05
11 3.68599613205767e-05
12 2.18429637716449e-05
13 1.29637467687612e-05
14 7.5934007464112e-06
15 1.01302942861322e-05
16 5.75126834442991e-06
17 3.25000357146837e-06
18 1.7876430746832e-06
19 9.61195281978375e-07
};
\addlegendentry{$\nu$}
\end{axis}

\end{tikzpicture}

%% file: figures/pgf/exact-error/street_solving_times.pgf
% This file was created with tikzplotlib v0.9.15.
\begin{tikzpicture}

\definecolor{color0}{rgb}{0.12156862745098,0.466666666666667,0.705882352941177}

\begin{axis}[
log basis y={10},
tick align=outside,
tick pos=left,
x grid style={white!69.0196078431373!black},
xlabel={Iteration},
xmajorgrids,
xmin=-0.95, xmax=19.95,
xtick style={color=black},
y grid style={white!69.0196078431373!black},
ylabel={Time (s)},
ymajorgrids,
ymin=0.0224883178643414, ymax=26.2269049830132,
ymode=log,
ytick style={color=black}
]
\addplot [semithick, color0, dashed, mark=*, mark size=1.75, mark options={solid}]
table {%
0 0.030999630689621
1 0.135999941267073
2 0.134000228717923
3 0.218999641947448
4 0.591000239364803
5 0.281999702565372
6 0.909999967552722
7 0.867999717593193
8 0.307000195607543
9 0.598999718204141
10 0.432999688200653
11 0.727000180631876
12 1.02399992756545
13 4.69200017396361
14 5.182000156492
15 2.46399978641421
16 4.84900055453181
17 8.69599990546703
18 10.899000428617
19 19.0260000759736
};
\end{axis}

\end{tikzpicture}

%% file: figures/pgf/exact-error/street_model_level_sets_sizes.pgf
% This file was created with tikzplotlib v0.9.15.
\begin{tikzpicture}

\definecolor{color0}{rgb}{0.12156862745098,0.466666666666667,0.705882352941177}
\definecolor{color1}{rgb}{1,0.498039215686275,0.0549019607843137}
\definecolor{color2}{rgb}{0.172549019607843,0.627450980392157,0.172549019607843}

\begin{axis}[
legend cell align={left},
legend style={
  fill opacity=0.8,
  draw opacity=1,
  text opacity=1,
  at={(0.91,0.5)},
  anchor=east,
  draw=white!80!black
},
tick align=outside,
tick pos=left,
x grid style={white!69.0196078431373!black},
xlabel={Iteration},
xmajorgrids,
xmin=-0.95, xmax=19.95,
xtick style={color=black},
y grid style={white!69.0196078431373!black},
ylabel={Proportion of pipes (\%)},
ymajorgrids,
ymin=-10, ymax=110,
ytick style={color=black}
]
\addplot [semithick, color0, dashed, mark=*, mark size=1.75, mark options={solid}]
table {%
0 0
1 0
2 0
3 0
4 0
5 0
6 0
7 0
8 0
9 0
10 0
11 0
12 0
13 0
14 0
15 0
16 0
17 0
18 0
19 0
};
\addlegendentry{M1}
\addplot [semithick, color1, dashed, mark=*, mark size=1.75, mark options={solid}]
table {%
0 0
1 9.87654320987654
2 22.2222222222222
3 34.5679012345679
4 46.2962962962963
5 46.2962962962963
6 56.7901234567901
7 66.0493827160494
8 73.4567901234568
9 79.0123456790123
10 79.0123456790123
11 82.7160493827161
12 86.4197530864197
13 89.5061728395062
14 91.9753086419753
15 91.9753086419753
16 93.8271604938272
17 95.0617283950617
18 96.2962962962963
19 97.5308641975309
};
\addlegendentry{M2}
\addplot [semithick, color2, dashed, mark=*, mark size=1.75, mark options={solid}]
table {%
0 100
1 90.1234567901235
2 77.7777777777778
3 65.4320987654321
4 53.7037037037037
5 53.7037037037037
6 43.2098765432099
7 33.9506172839506
8 26.5432098765432
9 20.9876543209877
10 20.9876543209877
11 17.2839506172839
12 13.5802469135802
13 10.4938271604938
14 8.02469135802469
15 8.02469135802469
16 6.17283950617284
17 4.93827160493827
18 3.7037037037037
19 2.46913580246914
};
\addlegendentry{M3}
\end{axis}

\end{tikzpicture}

%% file: conclusion.tex
\section{Conclusion}
\label{sec:conclusion}

In this paper, we set up a catalog of models for the hot water flow in
pipes of district heating networks.
For all entries of this catalog, we also derived suitable
discretizations to obtain finite-dimensional optimization problems for
the energy-efficient control of these networks that still ensures that
the demand of all customers are satisfied.
Based on these different models, we designed an iterative and adaptive
optimization method that automatically adapts the model level in the
catalog as well as the granularity of the discretization to
finally obtain a local optimal control that is feasible w.r.t.\ a
user-specified tolerance.
We show finite termination of this algorithm and present very
convincing numerical results that particularly show that we can now
solve realistic instances that are not solvable with state-of-the-art
commercial NLP solvers.

For our future work, we plan to extend our modeling and solution
approach to the case of instationary hot water flow modeling.
While we are confident that the overall ideas can be carried over to
this PDE-setting, this will most likely require some more technical
derivations compared to the ODE-case considered in this paper.

%%% Local Variables:
%%% mode: latex
%%% TeX-master: "dhn-adaptive"
%%% End:

%% file: acknowledgements.tex
\section*{Acknowledgments}
\label{sec:acknowledgements}

We acknowledge the support by the German
Bundesministerium für Bildung und Forschung within
the project \enquote{EiFer}.
Moreover, we are very grateful to all the colleagues within the EiFer
consortium for many fruitful discussions on the topics of this paper
and for providing the data.
We also thank the Deutsche Forschungsgemeinschaft for their
support within projects A05, B03, B08, and Z01 in the
Sonderforschungsbereich/Transregio 154 \enquote{Mathematical
  Modelling, Simulation and Optimization using the Example of Gas
  Networks}.

%%% Local Variables:
%%% mode: latex
%%% TeX-master: "dhn-adaptive"
%%% End:

%% file: dhn-adaptive.bib
@article{Mehrmann_et_al:2018,
    author = {Mehrmann, Volker and Schmidt, Martin and Stolwijk, Jeroen J.},
    year = {2018},
    title = {Model and Discretization Error Adaptivity Within Stationary Gas Transport Optimization},
    journal={Vietnam Journal of Mathematics},
    volume={46},
    number={4},
    pages={779-801},
    doi = {10.1007/s10013-018-0303-1},
}

@techreport{DomHLMMT21,
author={P. Domschke and B. Hiller and J. Lang and V. Mehrmann and R. Morandin and C. Tischendorf},
title={Gas Network Modeling: An Overview},
	note = {TRR 154 Preprint},
	year = {2021},
	url = {https://opus4.kobv.de/opus4-trr154},
}

@inproceedings{MehM19,
  author =	 {V. Mehrmann and R. Morandin},
  title =	 {Structure-preserving discretization for port-Hamiltonian descriptor systems},
  booktitle =	 {58th IEEE Conference on Decision and Control (CDC), 9.-12.12.19, Nice},
  pages =	 {6863--6868},
  publisher =	 {IEEE},
  year =	 2019,
  DOI = {10.1109/CDC40024.2019.9030180},
}

@article{Krug_et_al:2019,
  title     = {Nonlinear Optimization of District Heating Networks},
  author    = {Richard Krug and Volker Mehrmann and Martin Schmidt},
  journal   = {Optimization and Engineering},
  year      = {2021},
  volume    = {22},
  number    = {2},
  pages     = {783--819},
  doi       = {10.1007/s11081-020-09549-0},
}

@book{Stoer:2002,
    author = {Stoer, Josef and Bulirsch, Roland},
    year   = {2002},
    title  = {Introduction to Numerical Analysis},
    edition = {3},
    doi = {10.1007/978-0-387-21738-3},
    publisher = {Springer, New York, NY},
}

@incollection{Hauschild_et_al:2020,
  title        = {Port-{H}amiltonian modeling of district heating networks},
  author       = {Sarah-Alexa Hauschild and Nicole Marheineke and Volker Mehrmann and Jan Mohring and Arbi Moses Badlyan and Markus Rein and Martin Schmidt},
  year         = {2020},
  booktitle    = {Progress in Differential Algebraic Equations II},
  series       = {Differential-Albergaic Equations Forum},
  doi          = {10.1007/978-3-030-53905-4_11},
  editor       = {Reis, T. and Grundel, S. and Schöps, S.},
  publisher    = {Springer},
  preprint-url = {https://arxiv.org/abs/1908.11226},
  preprint-url = {https://opus4.kobv.de/opus4-trr154/frontdoor/index/index/docId/275},
}

@article{Hante_Schmidt:2019a,
  author   = {Falk M. Hante and Martin Schmidt},
  title    = {Complementarity-based nonlinear programming techniques for optimal mixing in gas networks},
  journal  = {EURO Journal on Computational Optimization},
  doi      = {10.1007/s13675-019-00112-w},
  volume   = {7},
  number   = {3},
  pages    = {299--323},
  year     = {2019},
  issn     = {2192-4414},
}

@article{Roland_Schmidt:2020,
  title    = {Mixed-Integer Nonlinear Optimization for District Heating Network Expansion},
  author   = {Marius Roland and Martin Schmidt},
  year     = {2020},
  journal  = {at - Automatisierungstechnik},
  note     = {Special Issue "Mathematical Innovations fostering the Energy Transition – Control, Optimization and Uncertainty Quantification"},
  doi      = {10.1515/auto-2020-0063},
  pubstate = {forthcoming},
  url-opto = {http://www.optimization-online.org/DB_HTML/2020/04/7752.html},
  url-opus = {https://opus4.kobv.de/opus4-trr154/frontdoor/index/index/docId/310},
}

@article{Stolwijk_Mehrmann:2018,
author = {Stolwijk, J. and Mehrmann, V.},
year = {2018},
month = {07},
pages = {231-266},
title = {Error Analysis and Model Adaptivity for Flows in Gas Networks},
volume = {26},
journal = {Analele Stiintifice ale Universitatii Ovidius Constanta, Seria Matematica},
doi = {10.2478/auom-2018-0027}
}

@incollection{NocSV09,
  title={Theory of adaptive finite element methods: an introduction},
  author={Nochetto, R.H. and Siebert, K.G. and Veeser, A.},
  booktitle={Multiscale, nonlinear and adaptive approximation},
  pages={409--542},
  year={2009},
  publisher={Springer},
  doi={10.1007/978-3-642-03413-8_12},
}

@Book{KonGMP03,
  Author =	 {M.~M. {Konstantinov} and D.~W. {Gu} and
                  V. {Mehrmann} and P.~H. {Petkov}},
  Title =	 {Perturbation Theory for Matrix Equations},
  ISBN =	 {0-444-51315-9/hbk},
  Pages =	 {xii + 429},
  year =	 2003,
  Publisher =	 {Amsterdam: North Holland}
}

@book{Ver13,
  title={A posteriori error estimation techniques for finite element methods},
  author={Verf{\"u}rth, R.},
  year={2013},
  publisher={OUP Oxford}
}

@inproceedings{Sandou_et_al:2005,
  author    = {{Sandou}, G. and {Font}, S. and {Tebbani}, S. and {Hiret}, A. and {Mondon}, C. and {Tebbani}, S. and {Hiret}, A. and {Mondon}, C.},
  booktitle = {Proceedings of the 44th IEEE Conference on Decision and Control},
  year      = {2005},
  doi       = {10.1109/CDC.2005.1583351},
  issn      = {0191-2216},
  pages     = {7372--7377},
  title     = {Predictive Control of a Complex District Heating Network},
}

@article{Verrilli_et_al:2017,
  author       = {{Verrilli}, F. and {Srinivasan}, S. and {Gambino}, G. and {Canelli}, M. and {Himanka}, M. and {Del Vecchio}, C. and {Sasso}, M. and {Glielmo}, L.},
  year         = {2017},
  journal      = {IEEE Transactions on Automation Science and Engineering},
  number       = {2},
  pages        = {547--557},
  title        = {Model Predictive Control-Based Optimal Operations of District Heating System With Thermal Energy Storage and Flexible Loads},
  volume       = {14},
  doi          = {10.1109/TASE.2016.2618948},
}

@article{Benonysson_et_al:1995,
  author       = {Benonysson, Atli and Bøhm, Benny and Ravn, Hans F.},
  year         = {1995},
  doi          = {10.1016/0196-8904(95)98895-T},
  issn         = {0196-8904},
  journal      = {Energy Conversion and Management},
  number       = {5},
  pages        = {297--314},
  title        = {Operational optimization in a district heating system},
  volume       = {36},
}

@article{Pirouti_et_al:2013,
  author       = {Pirouti, Marouf and Bagdanavicius, Audrius and Ekanayake, Janaka and Wu, Jianzhong and Jenkins, Nick},
  date         = {2013},
  doi          = {10.1016/j.energy.2013.01.065},
  issn         = {0360-5442},
  journaltitle = {Energy},
  pages        = {149--159},
  title        = {Energy consumption and economic analyses of a district heating network},
  volume       = {57},
}

@article{Rezaie_Rosen:2012,
  author       = {Rezaie, Behnaz and Rosen, Marc A.},
  date         = {2012},
  doi          = {10.1016/j.apenergy.2011.04.020},
  issn         = {0306-2619},
  journaltitle = {Applied Energy},
  pages        = {2--10},
  title        = {District heating and cooling: Review of technology and potential enhancements},
  volume       = {93},
}

@article{Schweiger_et_al:2017,
  author       = {Schweiger, Gerald and Larsson, Per-Ola and Magnusson, Fredrik and Lauenburg, Patrick and Velut, Stéphane},
  date         = {2017},
  doi          = {10.1016/j.energy.2017.05.115},
  issn         = {0360-5442},
  journaltitle = {Energy},
  pages        = {566--578},
  title        = {District heating and cooling systems – Framework for Modelica-based simulation and dynamic optimization},
  volume       = {137},
}

@article{Colella_et_al:2012,
  author       = {Colella, Francesco and Sciacovelli, Adriano and Verda, Vittorio},
  date         = {2012},
  doi          = {10.1016/j.energy.2012.03.043},
  issn         = {0360-5442},
  journaltitle = {Energy},
  note         = {The 24th International Conference on Efficiency, Cost, Optimization, Simulation and Environmental Impact of Energy, ECOS 2011},
  number       = {1},
  pages        = {397--406},
  title        = {Numerical analysis of a medium scale latent energy storage unit for district heating systems},
  volume       = {45},
}

@article{Verda_Colella:2011,
  author       = {Verda, Vittorio and Colella, Francesco},
  date         = {2011},
  doi          = {10.1016/j.energy.2011.04.015},
  issn         = {0360-5442},
  journaltitle = {Energy},
  number       = {7},
  pages        = {4278--4286},
  title        = {Primary energy savings through thermal storage in district heating networks},
  volume       = {36},
}

@article{Ben_Hassine_Eicker:2013,
  author       = {Hassine, Ilyes Ben and Eicker, Ursula},
  date         = {2013},
  doi          = {10.1016/j.applthermaleng.2011.12.037},
  issn         = {1359-4311},
  journaltitle = {Applied Thermal Engineering},
  note         = {Combined Special Issues: ECP 2011 and IMPRES 2010},
  number       = {2},
  pages        = {1437--1446},
  title        = {Impact of load structure variation and solar thermal energy integration on an existing district heating network},
  volume       = {50},
}

@article{Bordin_et_al:2016,
  author       = {Bordin, Chiara and Gordini, Angelo and Vigo, Daniele},
  date         = {2016},
  doi          = {10.1016/j.ejor.2015.12.049},
  issn         = {0377-2217},
  journaltitle = {European Journal of Operational Research},
  number       = {1},
  pages        = {296--307},
  title        = {An optimization approach for district heating strategic network design},
  volume       = {252},
}

@article{Dorfner_Hamacher:2014,
  author       = {Dorfner, J. and Hamacher, T.},
  date         = {2014},
  doi          = {10.1109/TSG.2013.2295856},
  issn         = {1949-3053},
  journaltitle = {IEEE Transactions on Smart Grid},
  number       = {4},
  pages        = {1884--1891},
  title        = {Large-Scale District Heating Network Optimization},
  volume       = {5},
}

@techreport{Borsche_et_al:2018,
  author = {Borsche, Raul and Eimer, Matthias and Siedow, Norbert},
  url    = {https://kluedo.ub.uni-kl.de/frontdoor/deliver/index/docId/5140/file/district_heating.pdf},
  date   = {2018},
  title  = {A local time stepping method for district heating networks},
}

@article{Rein_et_al:2018,
  author       = {Rein, Markus and Mohring, Jan and Damm, Tobias and Klar, Axel},
  date         = {2018},
  doi          = {10.1002/pamm.201800192},
  journaltitle = {PAMM},
  number       = {1},
  title        = {Parametric model order reduction for district heating networks},
  volume       = {18},
}

@techreport{Rein_et_al:2019,
  author       = {Rein, Markus and Mohring, Jan and Damm, Tobias and Klar, Axel},
  date         = {2019},
  url          = {https://arxiv.org/abs/1903.03342},
  title        = {Model order reduction of hyperbolic systems at the example of district heating networks},
}

@techreport{Rein_Mohring_et_al:2019,
  author        = {{Rein}, Markus and {Mohring}, Jan and {Damm}, Tobias and {Klar}, Axel},
  title         = {Optimal control of district heating networks using a reduced order model},
  url           = {http://publica.fraunhofer.de/documents/N-596673.html},
  date          = {2019}
}

@article{hart2011pyomo,
  author       = {Hart, William E and Watson, Jean-Paul and Woodruff, David L},
  publisher    = {Springer},
  year         = {2011},
  journal      = {Mathematical Programming Computation},
  number       = {3},
  pages        = {219--260},
  title        = {Pyomo: modeling and solving mathematical programs in {P}ython},
  volume       = {3},
  doi          = {10.1007/s12532-011-0026-8},
}

@book{hart2017pyomo,
  author    = {Hart, William E and Laird, Carl D and Watson, Jean-Paul and Woodruff, David L and Hackebeil, Gabriel A and Nicholson, Bethany L and Siirola, John D},
  publisher = {Springer},
  year      = {2017},
  title     = {Pyomo-optimization modeling in {P}ython},
  doi       = {10.1007/978-1-4614-3226-5},
}

@article{drud1994conopt,
  title={CONOPT—a large-scale GRG code},
  author={Drud, Arne Stolbjerg},
  journal={ORSA Journal on computing},
  volume={6},
  number={2},
  pages={207--216},
  year={1994},
  publisher={INFORMS}
}

@article{NUSSBAUMER2016496,
title = {Influence of system design on heat distribution costs in district heating},
journal = {Energy},
volume = {101},
pages = {496-505},
year = {2016},
issn = {0360-5442},
doi = {10.1016/j.energy.2016.02.062},
author = {T. Nussbaumer and S. Thalmann}
}

@book{Reid:1972,
title = {Riccati Differential Equations},
author = {William T. Reid},
series = {Mathematics in Science and Engineering},
publisher = {Elsevier},
volume = {86},
year = {1972},
issn = {0076-5392},
doi = {10.1016/S0076-5392(08)61166-2},
}

@article{Schmidt_et_al:2016,
  author       = {Martin Schmidt and Marc C. Steinbach and Bernhard M. Willert},
  title        = {High detail stationary optimization models for gas networks: validation and results},
  year         = {2016},
  journal      = {Optimization and Engineering},
  doi          = {10.1007/s11081-015-9300-3},
  issn         = {1389-4420},
  volume       = {17},
  number       = {2},
  pages        = {437--472},
  preprint-url = {http://www.optimization-online.org/DB_HTML/2014/10/4602.html},
}

@book{quarteroni2010numerical,
  title={Numerical Mathematics},
  author={Quarteroni, Alfio and Sacco, Riccardo and Saleri, Fausto},
  volume={37},
  year={2010},
  publisher={Springer Science \& Business Media},
  doi={10.1007/b98885},
}

@book{MR3443347,
  AUTHOR =	 {Bronshtein, I. N. and Semendyayev, K. A. and Musiol,
                  Gerhard and M\"{u}hlig, Heiner},
  TITLE =	 {Handbook of Mathematics},
  EDITION =	 {Sixth},
  PUBLISHER =	 {Springer, Heidelberg},
  YEAR =	 {2015},
  ISBN =	 {978-3-662-46220-1},
  MRCLASS =	 {00A20},
  MRNUMBER =	 {3443347},
  DOI     =      {10.1007/978-3-540-72122-2}
}
